# Ziele und Möglichkeiten des Einsatzes von Knowledge-Maps in mathematik-historischen Veranstaltungen

Nicola Oswald, Hannover; Sarah Khellaf, Hannover & Jana Peters, Hannover

*Zusammenfassung*: Die Vielschichtigkeit historischer Entwicklungen stellt bei der Vermittlung mathematik-historischer Inhalte eine Herausforderung dar. Um dieser zu begegnen, erprobten wir den Einsatz von Knowledge-Maps oder Wissenslandkarten. Der Fokus unseres Interesses liegt auf der Funktion der Maps zur Unterstützung des Meaningful Learning nach Ausubel (2000) sowie auf ihrer Rolle als Hilfsinstrument in mündlichen Prüfungen. In diesem Artikel wollen wir sowohl Praxiserfahrungen schildern, als auch Ergebnisse einer explorativen Studie in einem mathematik-historischen Seminar referieren.

*Abstract:* The multi-layered nature of societal developments poses a challenge in the teaching and learning of history of mathematics. As an attempt to tackle this challenge we experimented with the use of Knowledge-Maps. The focus of our interest lies on the use of these Maps to support meaningful learning (Ausubel, 2000) as well as on their role as facilitating tools in oral exams. In this article we seek to relay practical experiences we have made with Knowledge-Maps in our teaching and to report findings of an exploratory study we conducted in a seminar about history of mathematics.

## 1. Einleitung

Ausgangspunkt unserer Aktivitäten war der in vorangegangenen Semestern entstandene Eindruck, dass Studierende in Veranstaltungen mit geisteswissenschaftlichen Bezügen – und dies betrifft sowohl mathematik-historische, als auch mathematikdidaktische Inhalte – Schwierigkeiten bei der Organisation und Verknüpfung von erlerntem Wissen haben. Dies ist möglicherweise auf die Interdisziplinarität und die dadurch entstehende Vielschichtigkeit von Lehrangeboten zurückzuführen, die mathematische mit historischen- und sozialwissenschaftlichen Perspektiven verbinden. Bei der Vermittlung mathematik-historischer Inhalte stellen besonders die Prozesshaftigkeit und Verknüpfung historischer Entwicklungen untereinander eine Hürde dar.

In diesem Kontext hielten wir es für eine sinnvolle Unterstützungsmaßnahme, die Studierenden im Rahmen einer unbenoteten Studienleistung selbstständig Wissenslandkarten über die Inhalte von Vorlesungen und Seminaren erstellen zu lassen. Bei dieser Idee folgten wir unter anderem Resultaten von Barbro Grevholm

> […] it is clear that the maps give the student better opportunities to express her concept image. It obviously invites to more multidimensional answers than a sentence which in its form is linear. (2008, S. 6)

unter Beachtung der Hypothese von Nesbit und Adesope

> For instance, the mean effect sizes and confidence intervals across the subject categories suggest that concept mapping offers greater benefit in subject areas that are more saturated with verbal knowledge. However, the certainty of this interpretation is limited by significant heterogeneity and the small number of findings in the humanities, law, and social studies category. (2006, S. 428)

Im Folgenden werden diese Wissenslandkarten *Knowledge-Maps* oder kurz *Maps* genannt. Der gewählte Begriff soll die Nähe zum Concept- oder Mind-Mapping aufzeigen; unser Format verzichtet allerdings auf deren strenge gestalterische Vorgaben.

Zunächst möchten wir nun allgemein auf theoretische Grundlagen zur Vernetzung von Wissen eingehen. Hierbei unterscheiden wir zwischen einer mathematik-historischen und einer didaktischen Perspektive. Vor diesem Hintergrund entwickeln wir Forschungsfragen, welche wir im weiteren Verlauf dieses Beitrags näher beleuchten möchten.

### 1.1 Social Network Maps aus der mathematik-historischen Forschung

Social Network Analysis spielt seit einigen Jahren in den Geschichtswissenschaften eine entscheidende Rolle – so auch in der mathematik-historischen Forschung. Da dabei recht unterschiedliche Vorgehensweisen und sowohl qualitative als auch quantitative Methoden zum Einsatz kommen, plädieren Morten Reitmayer und Christian Marx in ihrem Übersichtskapitel „Netzwerkansätze in der Geschichtswissenschaft" im Grundlagenwerk *Handbuch Netzwerkforschung* (Stegbauer & Häußling, 2010) dafür, vielmehr von „Netzwerkansätzen zu sprechen" (Stegbauer & Häußling, 2010, S. 869). Damit wird bereits deutlich, dass die Erarbeitung und die Rekonstruktion von Netzwerken, sprich die Vernetzung von Interaktionen jeglicher Art, zunächst sehr frei gedacht werden kann. Zwei Wissenschaftler, die für die Verknüpfung von Geschichtswissenschaften und der so-zialen Netzwerkforschung argumentieren, sind etwa Marten Düring und, gerade im deutschsprachigen Raum, Martin Stark, die unter anderem das *Handbuch Historische Netzwerkforschung. Grundlagen und Anwendungen* (Düring, Eumann, Stark, & von Keyserlingk, 2016) herausgegeben haben und einen Überblick über aktuelle und vergangene Forschungsprojekte auf der Plattform *Historical Network Research* unter http://historicalnetworkresearch.org/ ermöglichen.

Grundsätzlich lassen sich hinsichtlich des Einsatzes von Netzwerken in der mathematik-historischen Forschung eine analysierende Funktion, etwa bei der quantitativen bibliometrischen Untersuchung großer Datensätze in Zitationsnetzwerken oder Kommunikationsnetzwerken[1], sowie eine illustrierende Funktion bei der Visualisierung von Netzwerken unterscheiden (z. B. Benstein & Oswald, 2016). Obwohl sowohl die Art der Umsetzung als auch die Auswirkungen auf den/die BetrachterIn kaum quantifizierbar und damit schwer zu bewerten sind, ist gerade die Veranschaulichung nicht zu unterschätzen:

> Perhaps the most central is the clear indication that imagery has, and has always had, a key role in network research. From the beginning images of networks have been used both to develop structural insights and to communicate those insights to others. (Freeman, 2000, S. 13)

Dieser Ansatz bildet Ausgangspunkt und Motivation für unseren Einsatz von Knowledge-Maps in mathematik-historischen Veranstaltungen. Im Folgenden möchten wir die Idee aus einer didaktischen Perspektive näher beleuchten.

## 1.2 Didaktische Forschung zu Map-Formaten in Bildungskontexten

Die Technik, Wissensbestände grafisch in sogenannten node-link Diagrammen zu ordnen, ist schon lange bekannt:

> Although the modern flourishing of node-link diagrams is often viewed as an entailment of Quillian's semantic networks (Quillian,1967), similar diagrams have been used for communication and learning since at least the 13th century (Sowa, 2000). (Nesbit & Adesope, 2006, S. 415)

In der zweiten Hälfte des 20. Jahrhunderts kam es unter anderem im Zuge der Entstehung kognitivistischer Lerntheorien zu einer Zunahme des Forschungsinteresses an dem Einsatz derartiger Visualisierungsverfahren in Bildungskontexten. Aktuelle Forschungsliteratur unterscheidet zwischen verschiedenen Arten von Mapping-Verfahren und zwischen den Theorieströmungen der Gestalt- und Kognitionspsychologie, auf die diese Verfahren begründet sein können. (Schukajlow & Leiß, 2012) In diesem Artikel folgen wir der Forschungstradition von Novak und Gowin (1984), die auf Basis von Ausubels kognitionspsychologischer Theorie des Meaningful Learning (Ausubel, 2000) sogenannte Concept Maps entwickelten und popularisierten. Zunächst soll nun eine Begriffsklärung erfolgen, bevor konkreter auf bisherige Forschungsergebnisse eingegangen wird.

Nesbit und Adesope (2006) verwenden in ihrer Meta-Analyse zum Lernen mit node-link Diagrammen den Begriff *Graphic Organizer*, um verschiedene Darstellungsformen von Wissen zu bezeichnen, die mit Ausubels Theorie in Verbindung gebracht werden. Hierunter fallen im Prinzip alle zweidimensionalen Visualisierungsverfahren[2], die durch geeignete Darstellungsmittel Beziehungen zwischen Begriffen oder Prozesse abbilden können. *Node-link Diagramme* sind eine Unterform, die – wie der Name nahelegt – auf Knoten und Verbindungslinien als primäre Darstellungsmittel zurückgreift. Diese werden wiederum unterteilt in *Concept Maps* nach Novak und Gowin (1984), *Knowledge Maps* nach Dansereau und Kollegen (siehe Nesbit & Adesope, 2006, S. 415) und weitere, anhand der jeweils geltenden Regeln bezüglich der erlaubten Darstellungsmittel. Zum Beispiel können Einschränkungen in Bezug auf die erlaubten Labels der Verbindungslinien oder in Bezug auf den Inhalt der Knoten gemacht werden. Des Weiteren können unterschiedliche Anforderungen in Bezug auf die Hierarchie der abgebildeten Informationen gestellt werden. Cañas et al. bezeichnen Concept Maps als semi-hierarchisch (2003, S. 13).

Wir arbeiten mit Graphic Organizers im weiteren Sinn. Der zentrale Unterschied zwischen unseren Knowledge-Maps und Concept Maps sowie den namensähnlichen Knowledge Maps nach Dansereau et al. besteht in der Abwesenheit strenger Gestaltungsregeln in unserem Arbeitsauftrag. Somit zählen auch Dokumente zu unseren Maps, die neben der standardmäßigen Beschriftung in den Knoten und an den Linien eine große Palette an verwendeten Darstellungsmitteln aufweisen. Unsere Entscheidung für diese sehr offene Definition von Knowledge-Map war eine Reaktion auf die hohe Komplexität der durch unsere Studierenden bearbeiteten Wissensbestände. Trotz der vordergründigen Formatunterschiede zwischen Concept Maps und den Produkten unserer Studierenden sind Forschungsergebnisse zur Wirkung dieses und ähnlicher Mapping-Verfahren aufgrund der gemeinsamen theoretischen Basis für unsere Aktivitäten relevant.

Im Folgenden werden nun einige für uns relevante Ergebnisse bisheriger Forschung referiert. In der von uns gesichteten Literatur finden sich neben Hypothesen zu möglichen Wirkmechanismen überwiegend quantitative Resultate (insbes. die Meta-Analyse von Nesbit und Adesope (2006)).

### 1.2.1 Bedingungen für Sinnvolles Lernen nach Ausubels Theorie

Ein Hauptziel des Einsatzes von Mapping-Verfahren in Bildungskontexten ist die Förderung Sinnvollen Lernens. In Ausubels Theorie ist der Anschluss an Vorwissen eine notwendige Bedingung für diese Art von Lernen. Insbesondere sollte auch auf sprachlicher Ebene ein Anschluss der neuen Inhalte an den aktuellen Stand des Lernenden möglich sein (Ausubel, 2000; Cañas et al., 2003, S. 6). Eine weitere notwendige Bedingung ist, dass die Lerner willentlich und aktiv Sinnvolles Lernen betreiben. (Ausubel, 2000; Cañas et al., 2003, S. 7)

### 1.2.2 Wofür werden Mapping-Verfahren in Bildungskontexten eingesetzt?

*Lernbegleitender Einsatz von Maps*

Lernbegleitend gibt es den Vorschlag, Concept Maps oder andere Graphical Organizers als *Advance Organizers* (Ausubel, 2000, S. 11) zu verwenden. Ein Advance Organizer wird bereits vor der Konfrontation mit dem zu lernenden Stoff eingesetzt, und soll den Zweck erfüllen, relevantes Vorwissen zu aktivieren, damit neue Inhalte in aktiver Auseinandersetzung mit dem Organizer daran angeknüpft werden können. In unserem Fall war die Idee, dass die Studierenden die in den einzelnen Veranstaltungssitzungen besprochenen Inhalte sukzessive und aufeinander bezogen in Form einer Knowledge-Map festhalten, sodass die Map die Funktion als Advance Organizer ab der zweiten Sitzung erfüllen kann.

Um Wissensspeicherung und -transfer anzuregen, scheint *Concept Mapping*, also das aktive Erstellen von Concept Maps, effektiver zu sein, als Texte zu lesen, Vorlesungen zu besuchen oder an Gruppendiskussionen teilzunehmen. Diese positive Wirkung konnte in einer Vielzahl von Kontexten und bei unterschiedlichen Lerngruppen nachgewiesen werden. Beim Lernen über konstruktive Aktivitäten mit Visualisierungsformaten wie Listen, Textausschnitten oder Outlines zeigten sich vergleichbare Effektstärken, weswegen möglicherweise die Konstruktivität der Aktivität, und nicht das spezielle Format von Concept Maps, den entscheidenden Faktor darstellt. (Nesbit & Adesope, 2006)

*Maps zur Zusammenfassung und zum Lernen (im Sinne von Behalten) von Information*

In der Literatur finden sich Vergleiche der Concept Mapping Technik mit anderen „summary formats" (Nesbit & Adesope, 2006), wo sie nach einer Lerneinheit zur Zusammenfassung und zum gezielten Lernen von Informationen genutzt wird. Eine Hypothese ist, dass die Übersichtlichkeit des Formats die Geschwindigkeit erhöhen könnte, mit der die Informationen (wiederholt) durchgesehen werden können.

Concept Maps wie auch unsere Knowledge-Maps sind in Bezug auf diesen Einsatz im Großen und Ganzen vergleichbar mit Formaten wie Listen oder eingerückten Outlines. Es wurden jedoch bei Concept Maps leicht stärkere positive Effekte auf das Lernen beobachtet, als bei den genannten Alternativen. (Nesbit & Adesope, 2006)

Ganz allgemein wird das Erstellen einer Concept Map von Novak und Gowin als eine kreative Handlung angesehen, die Reflexion über das eigene Wissen und Konzeptverständnis fördert:

> Novak and Gowin noted that the act of mapping is a creative activity, in which the learner must exert effort to clarify meanings, by identifying important concepts, relationships, and structure within a specified domain of knowledge. The activity fosters reflection on one's knowledge and understanding, providing a kind of feedback that helps students monitor their learning and, perhaps with assistance of teachers or peers, focus attention on learning needs. (Cañas et al., 2003, S. 22)

*Evaluation der Map: Bewertung*

Der Einsatz von Concept Maps als Bewertungsinstrument wird in der Literatur mit Skepsis betrachtet. Neben graphentheoretischen Verfahren, Maps zu scoren, zum Beispiel durch Zählen der Anzahl der verwendeten Knoten und Kanten, werden auch stärker inhaltliche Bewertungsverfahren vorgeschlagen, wobei die inhaltlichen den strukturellen Ansätzen in Bezug auf die Güte der Bewertung überlegen sind (Ley, 2014).

*Evaluation der Map: Diagnostik*

Der Einsatz zur formativen Evaluation bzw. Diagnostik, also zum Beispiel um kognitive Strukturen zu untersuchen oder Fehlvorstellungen zu identifizieren, wird hingegen von einigen Autoren empfohlen. In der Forschung wurde zum Beispiel die Kombination von Knowledge-Maps mit Interviews zur Aufdeckung kognitiver Strukturen der Interviewten genutzt (Grevholm, 2008). In pädagogischen Kontexten, vor allem in der Schule, finden solche Ansätze trotz bestehender Vorschläge (wie z.B. Kinchin, Hay & Adams, 2000) eher selten Anwendung (Ley, 2014).

Zwar geht es in der entsprechenden Literatur oft um die Nutzung von Mapping vor einer Lerneinheit, um Fehlvorstellungen und Vorwissen zu identifizieren (Cañas et al., 2003, S. 23), auf theoretischer Ebene spricht aber erst einmal nichts gegen einen Einsatz zu beliebigen Zeitpunkten im Lernprozess, um kognitive Strukturen abzubilden.

*Einsatz als Kommunikationshilfe*

Vorstrukturierte Maps können als Kommunikationsstütze für sprachlich schwächere Lerner dienen; dieser Effekt war bei sprachlich starken Lernern nicht feststellbar (Nesbit und Adesope, 2006). Des Weiteren kann bei Gruppenaktivitäten die Map als Fokuspunkt der Gruppenkommunikation hilfreich sein. In ihrer Meta-Studie können Nesbit und Adesope (2006) zwar keine Effekte nachweisen, untermauern diese Hypothese jedoch anhand inhaltlicher Überlegungen:

> We believe that the potential advantages of concept maps in cooperative learning are sensitive to both the nature of the task and the training of participants in cooperative methods. Tasks that require rapid communication of complex, non-hierarchically structured information are more likely to benefit from the conceptually integrated propositions represented in the concept map format. (Nesbit und Adesope, 2006, 431)

**1.3 Anliegen des Artikels**

Auf Basis der bisher referierten Informationen und Gedanken sollen im Verlauf des restlichen Artikels zum einen Praxiserfahrungen beschrieben werden, die wir bisher mit dem Einsatz von Knowledge-Maps gemacht haben. Zum anderen wollen wir die Ergebnisse einer explorativen Studie referieren, die wir in einem mathematik-historischen Seminar am Standort Hannover durchgeführt haben. Dabei sollen folgende Forschungsfragen näher beleuchtet werden:

- Welche Lernhandlungen regte der im Seminar gestellte Knowledge-Map-Arbeitsauftrag an?

- Welche in der auf Ausubels Theorie basierenden Literatur genannten Wirkungen der Arbeit mit Graphic Organizers – nachgewiesene oder antizipierte – können wir in unseren Daten wiederfinden?

- Was für Erfahrungen haben wir mit dem Einsatz von Knowledge-Maps als Unterstützung für mündliche Prüfungen gemacht? Was kann das Dokument Knowledge-Map in mündlichen Prüfungen leisten?

**2. Lehrveranstaltungen**

Unsere bisherigen Erfahrungen mit Knowledge-Maps entstanden in Veranstaltungen zur Entwicklung mathematischer Denkweisen vor gesellschaftshistorischem Hintergrund, welche an der Bergischen Universität Wuppertal sowie an der Leibniz Universität Hannover abgehalten wurden. Im Folgenden werden in Abschnitt 2.1 zunächst Inhalte und didaktische Hürden von Vorläuferveranstaltungen aus Wuppertal skizziert, bevor die aus diesen Erfahrungen entstandenen Ziele und Ideen des Einsatzes von Knowledge-Maps formuliert werden. Anschließend gehen wir in Abschnitt 2.2 auf das Hannoveraner Seminar aus dem Wintersemester 2017/18 ein, welches die Daten für die anschließende Studie lieferte.

**2.1 Vorlesungen in Wuppertal**

Als Vorläufer und Kontext der anschließend beforschten Veranstaltung soll kurz die Vorlesung „Entwicklung mathematischer Denkweisen" (4 SWS) vorgestellt werden, welche in den Sommersemestern 2016 und 2017 an der Bergischen Universität Wuppertal durchgeführt wurde. Diese Veranstaltung bildet, alternierend mit der Vorlesung „Ausgewählte Kapitel der Mathematikgeschichte", in Wuppertal einen festen Bestandteil des mathematik-historischen Pflichtmoduls für Studierende im Master of Education (Lehramt Gymnasium). Das Modul wird mit einem ergänzenden mathematischen oder mathematik-historischen Seminar abgeschlossen, wobei sich die numerische Bewertung des gesamten Moduls aus Vortrag und Hausarbeit im Rahmen dieses Seminars ergibt. Die Tatsache, dass die Anwesenheit der Studierenden in der Vorlesung nicht verpflichtend gefordert werden soll, stellt Dozierende vor eine Herausforderung: Wie kann eine adäquate unbenotete Studienleistung erhoben werden und zugleich das Interesse und die Motivation der Studierenden für Mathematikgeschichte angeregt und gefördert werden?

Ausgehend von dem in Abschnitt 1.1 beschriebenen Vorwissen zu Social Network Maps in der mathematik-historischen Forschung erprobten wir im Sommersemester 2016 zur Bewältigung dieser Herausforderung eine auf Knowledge-Maps basierende Methode, welche in den darauffolgenden Semestern, vor allem in Hinblick auf die Einbeziehung der lerntheoretischen Überlegungen aus Abschnitt 1.2, weiterentwickelt werden konnte.

2.1.1 Inhalte der Vorlesung

Zur Illustration unserer Herangehensweise sowie der ersten Wuppertaler Erfahrungen, gehen wir zunächst auf die Inhalte und Lernziele der Vorlesungen ein. Diese ähnelten sich 2016 und 2017 bis auf wenige Details bezüglich der Schwerpunktsetzung. Als Leitfaden für beide Vorlesungen wurden folgende Fragen kritisch behandelt: Was genau kann unter einer ‚mathematischen Denkweise' verstanden werden und inwiefern ist es überhaupt möglich, eine solche Denkweise zu charakterisieren? Wann kann in der historischen Betrachtung mathematischer Denkweisen von einem Entwicklungsschritt oder einer Weiterentwicklung gesprochen werden? Den Anliegen der Lehramtsstudierenden entgegenkommend, wurde zu-sätzlich an wiederholter Stelle darauf eingegangen, ob und wie behandelte mathematik-historische Sachverhalte in den heutigen Mathematikunterricht einbezogen werden könnten bzw. inwieweit deren Einbeziehung didaktischen Prinzipien entsprechen kann.

Um einen Eindruck vom behandelten stofflichen Umfang zu geben, ohne inhaltliche Details zu vertiefen, geben wir hier einen stichpunktartigen Überblick über die einzelnen Themengebiete, welche weitgehend in chronologischer Reihenfolge behandelt wurden: Nach einer einführenden Sektion zur frühen Entwicklung von Individualzeichen, der mathematischen Entwicklung in Mesopotamien anhand von ausgewählten Tontafeln, Persönlichkeiten und ihnen zugeschriebene Errungenschaften aus der griechischen und chinesischen Antike, folgte ein zeitlicher und geografischer Sprung ins nach-neuzeitliche Europa. Anhand von Rechenweisen nach Adam Ries sowie Resultaten von Pierre de Fermat, Nicole-Reine Lepaute, Carl Friedrich Gauß und Sophie Germain reflektierten wir Weiter- und Neuentwicklungen der Arithmetik und Zahlentheorie. Mit der Frage nach der Ausbildung von Disziplinen innerhalb der Mathematik beschlossen wir die Vorlesung durch einen Blick auf die Entwicklung der „Geometrie der Zahlen", stark geprägt durch Hermann Minkowski im ausgehenden 19. Jahrhundert.

Insgesamt wurde folglich eine zeitlich recht weite Spanne, naturgemäß mit Brüchen, abgedeckt, wobei die einzelnen Themen immer wieder auf Einzelpersonen bezogen sowie an meist arithmetischen bzw. zahlentheoretischen Beispielen auf mathematischer Ebene erläutert wurden. Um einer punktuellen Hervorhebung einzelner historischer Ereignisse entgegen zu wirken, lag zudem ein Hauptaugenmerk auf der Vermittlung eines kritischen Verständnisses für mathematik-historische Geschichtsschreibung. Unter anderem sollten Entwicklungen vor dem jeweiligen gesellschaftlichen und politischen Hintergrund reflektiert werden.

### 2.1.2 Lernziele und Vorgaben für die Knowledge-Maps

Der letztgenannte Anspruch spiegelt sich in den Lernzielen der Vorlesung wider: Die Studierenden sollten lernen, mathematische Entwicklungen vor gesellschaftlich-historischem Hintergrund zu reflektieren, Zusammenhänge zwischen den Inhalten der Vorlesung zu erkennen und sich eigenständig anhand der behandelten Themen einen chronologischen Überblick zu erarbeiten. Diese Ziele waren zu Beginn des Semesters implizit durch Hinweise zur individuellen Erstellung einer Knowledge-Map (siehe Abbildung 1), gegeben. Auf diesem Merkblatt wurde den Studierenden zudem erläutert, dass deren Abgabe sowie die Beantwortung eines Fragebogens am Ende des Semesters unbewertet, jedoch verpflichtend sei, und es wurden weiterführende allgemeine Informationen zu dem Visualisierungswerkzeug Concept-Map angeboten. Die Gestaltungsmöglichkeiten wurden allerdings, entgegen den klassischen Vorgaben bei der Erstellung einer Concept-Map, bewusst freigehalten. Diese Freiheit wurde von den Studierenden sehr vielfältig genutzt. Die schließlich abgegebenen Maps wiesen in beiden Semestern nicht nur sehr unterschiedliche Herangehensweisen hinsichtlich der Gestaltung auf, sondern zeigten in großen Teilen auch eine zunächst unerwartete, da es sich um eine unbenotete Studienleistung handelte, Komplexität bezüglich der Detailgenauigkeit sowie der Verknüpfungen von einzelnen Inhalten (siehe Abbildung 2). Am Ende des Semesters wurde von allen Studierenden der erwähnte Fragebogen zu Zusammenhängen der Vorlesungsinhalte ausgefüllt. Dieser war bewusst frei formuliert, konnte jeweils einer Map zugeordnet werden und sollte der Dozentin einen Eindruck Wissensstand der Studierenden vermitteln. Beispielsweise wurde die Frage „Fermat und Mersenne sind Teil einer neuen Form eines Austausches von Mathematik An welche Aspekte erinnern Sie sich?" und „Warum wird Gauß weithin als Begründer der modernen Zahlentheorie betrachtet?" gestellt.

Abb. 1 Richtlinien zur Erstellung der Knowledge-Maps, in den SoSe 2016 und 2017 in Wuppertal.

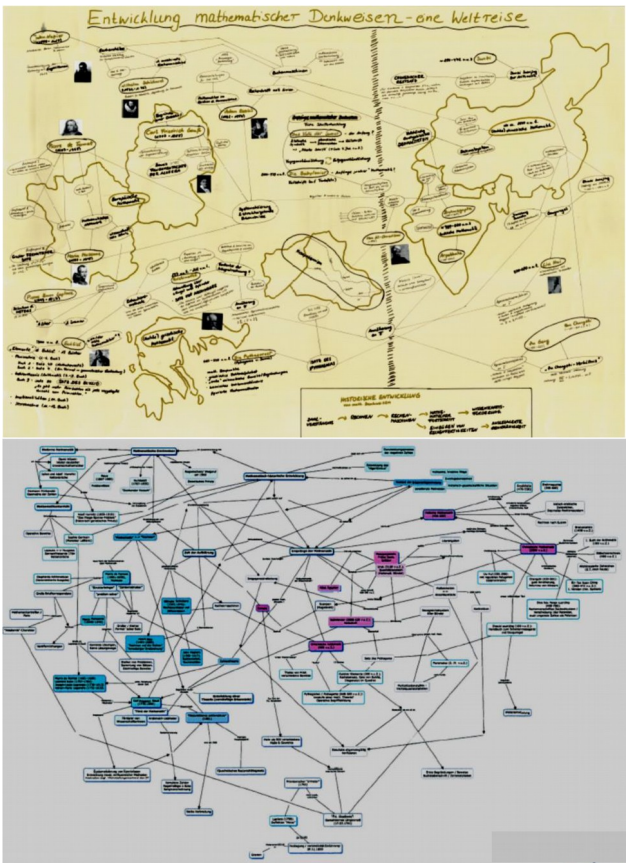

Abb. 1 Beispielkarten aus den SoSe 2016 und SoSe 2017 in Wuppertal.

### 2.1.3 Idee einer möglichen Nutzung der Maps zur Leistungsbewertung

Im Rahmen der Masterarbeit (Schumacher, 2017) wurde im Sommersemester 2017 ein Versuch unternommen, ein Bewertungsschema für Knowledge-Maps zu entwerfen, das zum Zwecke der Leistungsbewertung eingesetzt werden könnte. Bei ihrer Herangehensweise orientierte sich die Studentin insbesondere an Veröffentlichungen des AK Vernetzung im Mathematikunterricht (Brandl M., Brinkmann A., Maaß J. & Siller H.-S., 2016), und untersuchte die entstandenen Knowledge-Maps aus beiden Semestern hinsichtlich der Gütekriterien Funktionalität, Kognition, Technik und Gestaltung.

Die Bewertung erfolgte anhand eines von der Studentin entwickelten Punktesystems. An dieser Stelle muss betont werden, dass insbesondere die kreative Leistungsbewertung (Gestaltung und Technik) eine Herausforderung darstellte und stark intuitiv geprägt war. Dementsprechend können die Ergebnisse lediglich als Tendenzen betrachtet werden und sollen hier nicht weiter vertieft werden. Jedoch möchten wir kurz auf Hinweise zu Zusammenhängen zwischen der Qualität der Fragebögen und der Knowledge-Maps eingehen, welche von der Studentin folgendermaßen untersucht wurden:

> […] im Zusammenhang der „Kognition" und in Bezugnahme auf die „Funktionalität" [sollte] die Vollständigkeit der Maps betrachtet werden, indem Antworten auf die Fragen des Fragebogens auf den Landkarten zu finden versucht werden. (Schumacher, 2017, S. 36)

Insgesamt zeigte sich, dass Studierende, deren Knowledge-Maps hinsichtlich Funktionalität und Kognition eine hohe Punktzahl erzielen konnten auch eine gute bis sehr gute Leistungsbewertung erlangten. Die Umkehrannahme, dass als unvollständig bewertete Knowledge-Maps Rückschlüsse auf eine geringere Leistung hinsichtlich der Fragebögen zuließe, traf jedoch nicht zu. Allerdings konnten Studierenden mit geringer Leistungsbewertung der Fragebögen tendenziell gering bewertete Knowledge-Maps zugeordnet werden.

### 2.1.4 Wuppertal als Vorläufer und Motivation für Hannover

Aus Dozentinnensicht konnte nach den beiden Vorlesungen in Wuppertal festgehalten werden, dass die Abgabe von Knowledge-Maps zur Erhebung einer unbenoteten Studienleistung aus didaktischer Sicht durchaus vielversprechend war. Es konnten Tendenzen dafür erkannt werden, dass die Erstellung der Map den Studierenden als Hilfestellung bei der Strukturierung und der Verknüpfung mathematikhistorischer Inhalte diente. Dass die Knowledge-Maps in großen Teilen eine hohe Detailtreue und komplexe Verzweigungen der Inhalte aufwiesen, interpretierten wir als Anzeichen einer hohen Motivation der Studierenden. Diese Beobachtungen bildeten die Ausgangslage für die im Folgenden beschriebene Weiterentwicklung des gezielten Einsatzes von Knowledge-Maps in mathematik-historischen Veranstaltungen.

## 2.2 Das Seminar in Hannover

Im Wintersemester 2017/18 wurde an der Leibniz Universität Hannover das Seminar „Entwicklung mathematischer Denkweisen im gesellschafts-historischen Kontext" (2 SWS) gehalten. Die Inhalte des Seminars wurden in starker Anlehnung an die zuvor beschriebenen Vorlesungen gewählt. Soweit möglich waren die einzelnen Themen sowohl chronologisch als auch lokal geographisch geordnet: Von babylonischer Mathematik in Mesopotamien, Schwerpunkten der griechischen und chinesischen Antike, Beispielen islamischer Mathematik bis ins Mittelalter hin zu ausgewählten europäischen MathematikerInnen und ihren Errungenschaften ab dem 16. Jahrhundert (insbesondere Adam Ries, Pierre de Fermat, Carl Friedrich Gauß und Sophie Germain). Neben der Darstellung fachmathematischer Sachverhalte lag ein Hauptaugenmerk auf der Reflexion der Themen vor dem jeweiligen gesellschaftlichen, politi-schen, beziehungsweise kulturellen Hintergrund, beispielsweise die Arbeit von Sophie Germain vor dem Hintergrund der gesellschaftlichen Stellung der Frau.

Hinsichtlich der Art der Aufarbeitung sowie der individuellen Schwerpunktsetzung innerhalb der einzelnen Themenbereiche hatten die Studierenden der Natur der Veranstaltungsart gemäß einen großen Gestaltungsspielraum. So entschieden sich beispielsweise die Studierenden, die das Thema „Pierre de Fermat: Akademie der Briefe rund um Marin Mersenne?" vorstellten, für einen starken Fokus auf den so genannten Letzten Satz von Fermat. Sie schlugen einen unerwartet weiten Bogen bis in das Jahr 1994 und den endgültigen Beweis durch Andrew Wiles. Ähnlich verhielt es sich mit der Referatsgruppe zum Thema „Historische Rechenmaschinen". Dieses war ursprünglich von der Dozentin chronologisch im Hinblick auf mathematische Denkweisen als „Auslagern von Grundrechenarten" auf Maschinen, hauptsächlich im 17. und 18. Jahrhundert, vorgesehen. Die Studierenden entschieden sich jedoch für eine umfassende Darstellung der Geschichte der Rechenmaschinen anhand von umsichtig gewählten Beispielen, ausgehend von antiken Rechentafeln bis hin zur Entwicklung des ZUSE Computers.

2.2.1 Lernziele und Erwartungshorizont

Damit erfüllten diese Studierenden bereits ein prinzipielles Ziel von Seminaren: Durch die selbstgewählte Herangehensweise an ihre Themen bewiesen sie Autonomie und Eigeninitiative. Diese Aspekte spielten auch explizit bei den speziell für dieses Seminar geforderten Lernzielen eine Rolle, welche in der Veranstaltung formuliert wurden: Die Studierenden sollten ein Verständnis für historische Entwicklungsschritte sowie gesellschaftliche Einflüsse auf mathematische Denkweisen entwickeln, Diskussionen der Inhalte und Forschungsansätze im Anschluss an die Seminarbeiträge anregen und führen, sensibilisiert werden für den Umgang mit historischen Themen und Historiographie und deren Verwendung, und darüber hinaus mathematische Sachverhalte aufarbeiten, insbesondere hinsichtlich der eigenen Referate, sowie die kritische Auseinandersetzung mit historischen Quellen und Bezügen zeigen. Diese Anforderungen wurden sowohl schriftlich im Rahmen einer Einführungsveranstaltung zu Beginn des Semesters als auch verbal in Diskussionen mit den Studierenden im Laufe der Veranstaltung benannt. Die Kommunikation des Erwartungshorizonts wurde bewusst so transparent und nachvollziehbar wie möglich gestaltet.

2.2.2 Studien- und Prüfungsleistung

Hinsichtlich der im Seminar zu erbringenden Leistungen hatten die Knowledge-Maps eine Doppelfunktion. Zum einen fungierten sie, neben den thematisch fokussierten Referaten, als *unbenotete Studienleistungen*. Der Auftrag lautete, die Knowledge-Maps seminarbegleitend zu erstellen (siehe Abb. 3). Außerdem sollte in der Map der gesamte Stoff der Veranstaltung abgedeckt sein. Intendiert wurde dadurch die Nutzung der Maps als lernbegleitende Vernetzungs- und Reflexionshilfe (vgl. Abschnitt 1.2). In diesem Kontext gingen wir davon aus, dass gemäß Ausubels Voraussetzung relevantes Vorwissen auf Seiten der Studierenden vorhanden war. Es wurden in den Seminarsitzungen kaum historische Originaltexte gelesen, die eine sprachliche oder inhaltliche Hürde hätten darstellen können, und zu den behandelten Konzepten sollten Studierende i.A. Vorwissen besitzen. Zum anderen sollten die Knowledge-Maps im Sinne einer Gesprächsgrundlage als *Hilfsmittel für die benoteten mündlichen Prüfungen* am Ende des Semesters dienen. Als Ergebnisse eines kreativen Schöpfungsprozesses erwarteten wir, dass die Maps der Studierenden Gesprächsanlässe schaffen würden. Darüber hinaus waren wir an ihrem Potential als Übersicht schaffendes Kommunikationsmittel interessiert.

Abb. 3: Hinweise zu Knowledge-Maps, Einführungsveranstaltung im Seminar, WS 2017/18

Bereits zu Beginn des Semesters wurden die Studierenden auf diese beiden Rollen der Knowledge-Map hingewiesen. Begleitend wurden ihnen, ähnlich wie in den Vorlesungen in Wuppertal, in der Einführungsveranstaltung einige grundlegende Hinweise zur Erstellung von Knowledge-Maps (siehe Abbildung 3) sowie Beispiele aus den vorangegangenen Veranstaltungen in Wuppertal zur Verfügung gestellt (siehe Abbildung 2). Die Knowledge-Maps sollten bis Ende Januar 2018 abgegeben werden. Die mündlichen Prüfungen fanden am 06., 19. und 20. Februar 2018 statt.

Eine dritte Funktion hatten die Knowledge-Maps hinsichtlich der *Vorbereitung der mündlichen Individualprüfungen seitens der Dozentin*. Die Gesamtprüfung bestand aus einem Block mündlicher Individualprüfungen, die in einem maximal 15-minütigen Gespräch zu themenübergreifenden Inhalten des Seminars bestanden, und einem zweiten Prüfungsteil – einer mündlichen Gruppenprüfung der Mitglieder einer Referatsgruppe zu deren Vortragsthema. Das im Weiteren beschriebene Vorbereitungsverfahren bezog sich lediglich auf die vorgeschalteten Individualprüfungen, nicht auf das Gruppengespräch.

Die Vorbereitung der Individualprüfungen bestand zunächst aus der Festlegung des Erwartungshorizontes (Abbildung 4), der sich an den Lernzielen des Seminars (Abschnitt 2.2) orientierte. Darin wurden vier Anforderungskategorien festgelegt. Die Notenstufe eines Prüfungsergebnisses sollte der Anzahl der durch den/die zu Prüfenden nicht erfüllten Anforderungen plus eins entsprechen. In einem zweiten Schritt konnte der Grad der Erfüllung der Anforderungen durch die Wahl von Notenabstufungen aus-

**(allgemeiner) Erwartungshorizont an die Studierenden**

1. Die Studierenden sollen zeigen, dass sie mathematische Inhalte vor kulturhistorischem / sozialem Hintergrund reflektieren können.
2. Die Studierenden sollen inhaltliche Zusammenhänge zwischen den einzelnen Beiträgen des Seminars herstellen können (vor dem Kontext der Entwicklung mathematischer Denkweisen).
3. Die Studierenden sollen zeigen, dass sie sich mit (mathematik-historischen) Methoden kritisch auseinander gesetzt haben. (Dazu gehört auch, Fehlvorstellungen zu reflektieren).
4. Die Studierenden sollen zeigen, dass sie mit der Mathematik der Gebiete ihres Seminarbeitrags vertraut sind und sicher umgehen können.

Abb. 4 Erwartungshorizont der mündlichen Prüfung

gedrückt werden.

Als Hilfe zur konkreten Prüfungsdurchführung wurden ferner standardisierte Fragen zu den Themen, Gestaltung und Inhalt entwickelt. Die ursprüngliche Fragensammlung findet sich in Abbildung 5. Sie wurde aber nach Sichtung der ersten Maps erweitert durch die Frage nach dem Vorhandensein der Darstellung einer chronologischen Entwicklung, die der Kategorie Inhalt zugeordnet wurde. Die Fragen zum

1. **Gestaltung**
   1. Erklären Sie kurz Ihre grundsätzliche Vorgehensweise bei der Erstellung der Map?
   2. Haben Sie für Ihre Map bewusst ein Zentrum gewählt?
      a) Wenn ja, warum?
      b) Wenn nein, warum nicht?
   3. Welche visuellen Mittel (Farben, Bilder, Pfeife, ect.) haben Sie wofür verwendet? Erläutern Sie kurz warum?
2. **Inhalt (an der jeweiligen Map orientiert)**
   2.1. Erklären Sie den von Ihnen beschriebenen Zusammenhang Y zwischen den Themen X und Z?
   2.2. Warum haben Sie Thema X nicht in Ihre Map aufgenommen? Was wurde zu diesem Thema besprochen? Wo könnte dieses in Ihrer Map eingeordnet werden?
   2.3. allgemein: Nennen Sie Entwicklungsschritte mit Bezug zu ….

Abb. 5 Prüfungsfragen

Inhalt wurden von der Dozentin dafür genutzt, im Vorfeld einer Prüfung die jeweilige Map einmal zu sichten, und weitere, individuell auf diese Map angepasste Fragen abzuleiten.

Im Prüfungsgespräch sollte zunächst anhand von Fragen zur Gestaltung eine Begründung zum Aufbau der Map sowie zu speziell gewählten gestalterischen Elementen der Map-Struktur gegeben werden, welche Rückschlüsse auf die erarbeitete Wissensstruktur des/der jeweiligen Studierenden bezüglich des gesamten behandelten Stoffes im Seminar liefern sollte. Mit den individuellen Fragen zum Inhalt sollten anschließend ausgewählte Verbindungen fokussiert werden, welche auf der gegebenen Knowledge-Map zu erkennen oder explizit nicht zu erkennen waren. Hier galt es, herauszufinden inwieweit die Studierenden Verknüpfungen zwischen mathematik-historischen Entwicklungsschritten herstellen konnten bzw. warum sie diese nicht dargestellt hatten.

Dieses bisher eher knapp beschriebene Vorgehen der Prüfungsvorbereitung durch die Dozentin soll nun anhand eines Beispiels – hier und im Folgenden KL genannt – illustriert werden. In Abbildung 6 ist die entsprechende Knowledge-Map abgebildet. Gemäß obigen Verfahrens unterzog die Dozentin die Knowledge Map KL ca. einen Tag vor der zugehörigen mündlichen Prüfung einer Voranalyse, indem sie die Map auf die vorbereiteten allgemeinen Fragen untersuchte.

In unserem Beispiel fiel neben der insgesamt sehr kleinteiligen und detailreichen Gestaltung insbesondere die selbstständig entwickelte, strukturgebende Einteilung gesellschaftlicher Entwicklungsstufen auf. Den „Entwicklungen in der Mathematik" (Farbe: blau) wird eine Liste „Gesellschaft bedarf Mathematik für …" gegenübergestellt, welche mithilfe farblicher Zuordnung die Themen des Seminars strukturiert. Beispielsweise wird der frühen „Ägyptischen Mathematik" die benötigte „Zentralstaatkoordinierung" zugeschrieben. Zusätzlich wird eine Vielzahl von Verbindungen zwischen einzelnen Themen gezogen. So wird etwa eine „Weiterentwicklung von …"-Beziehung zwischen der griechischen Antike (Dreiecks- und Pentagonzahlen) und der islamischen Mathematik (bildlichen Algebra) hergestellt. Zu beiden Seminarbeiträgen wurde explizit die Rolle von figurierter Mathematik diskutiert.

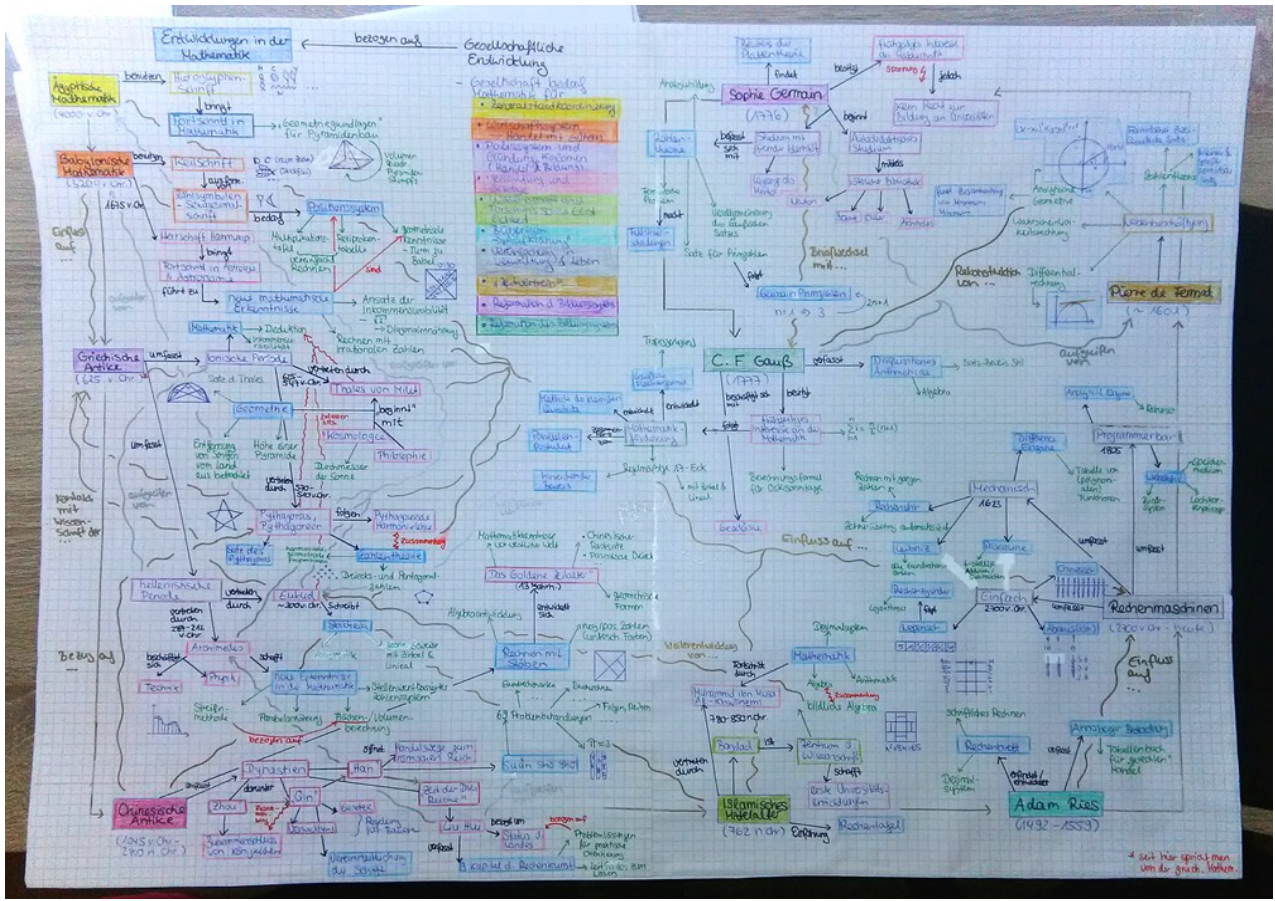

| Knowledge-Map KL 06.02.2018 | kein Zentrum, klassisches Mapping mit vielen Verzweigungen extra: Liste zu Kriterien gesellschaftlicher Entwicklungen, chronologisch? | Verbindung Y: Figurierte Mathematik zwischen X: Pythagoreer Z: bildliche Algebra (islamische Mathematik) | aufwendige kleinteilige Gestaltung, Struktur durch Farbigkeit und Verknüpfungen, Ordnung durch Liste zu gesellschaftlichen Entwicklungsschritten |
|---|---|---|---|

Abb. 6 Map zu Fallstudie II und Auszug Prüfungsvorbereitung der Dozentin

Die genannten Beobachtungen der Dozentin wurden als strukturierte Hilfsnotizen schriftlich festgehalten (siehe Abbildung 6 unten) und fungierten während der Prüfung als Gesprächsleitfaden.

Nach Abschluss der Gesamtprüfung aber vor Bekanntgabe der Noten wurden die Studierenden gebeten, einen kleinen Fragebogen auszufüllen (siehe Abbildung 7). Dieser war bewusst niedrigschwellig und kurz konzipiert und sollte uns einen Eindruck vom Umgang der Studierenden mit den Knowledge-Maps vermitteln.

Abb. 7 Beispiel eines Fragebogens

## 3. Erhebungen und Analyse

Insgesamt haben 25 Studierende am Seminar „Entwicklung mathematischer Denkweisen im gesellschafts-historischen Kontext" im Wintersemester 2017/18 teilgenommen. Davon haben 21 Studierende die mündliche Prüfung gemacht. Unsere Analyse beschränkt sich im Folgenden auf Datenmaterial, das durch diese 21 Studierenden generiert wurde. Das Material besteht aus den abgegebenen Knowledge-Maps, Audioaufzeichnungen der mündlichen Prü-fungen und den von den Studierenden im Anschluss zur Prüfung ausgefüllten Fragebögen.

### 3.1 Daten

Die von den Studierenden abgegebenen **Knowledge-Maps** wurden vor der Analyse digitalisiert.

Vom **ersten Teil der mündlichen Prüfungen**, der sich auf die Knowledge-Maps bezog, wurde mit Einverständnis der Studierenden eine Audioaufnahme gemacht, die anschließend transkribiert wurde. In unserer explorativen Analyse haben wir im Wesentlichen die Antworten der Studierenden auf die erste Frage „Was haben Sie sich bei der Erstellung der Knowledge-Map gedacht, wie sind Sie vorgegangen?" berücksichtigt. Darüber hinaus wurde die Interaktion der Studierenden mit der Map (auf die Map schauen und deuten) während dieses gesamten ersten Prüfungsteils qualitativ (sehr viel, viel, wenig) durch eine Prüfungsbeisitzende protokolliert.

Der **Fragebogen**, den die Studierenden nach der Prüfung ausfüllten, sollte einen Eindruck vom Umgang der Studierenden mit den Knowledge-Maps vermitteln. Die Fragen widmen sich dem Zeitpunkt der Erstellung der jeweiligen Map, dem zeitlichen Aufwand sowie der individuellen begründeten Einschätzung der Sinnhaftigkeit der Aufgabe, parallel zum Seminar eine Knowledge-Map zu erstellen. Die beantworteten Bögen wurden ebenfalls digitalisiert.

### 3.2 Analysemethode

Zur Untersuchung unserer Forschungsfragen orientierten wir uns an der computergestützten qualitativen Inhaltsanalyse nach Mayring (2008) und Kuckartz (2016)[3].

Die **Knowledge-Maps** codierten wir nach einer inhaltlichen und einer gestalterischen Dimension. Bei der Codierung hinsichtlich der inhaltlichen Dimension wurden nur Aspekte codiert, die die Struktur der Knowledge-Map maßgeblich bestimmt haben. Wurde beispielsweise der Aspekt Zeit genutzt, um die Map chronologisch zu strukturieren, wurde ein entsprechender Code zugewiesen. Wurden hingegen Zeitangaben ohne Einfluss auf die Struktur angebracht wurde Zeit nicht codiert. Hinsichtlich der gestalterischen Dimension wurden strukturgenerierende Designelemente codiert.

Die Codes für beide Dimensionen wurden sowohl deduktiv als auch induktiv aus den Knowledge-Maps entwickelt. Dabei bezogen wir uns bei der Charakterisierung der deduktiven Codes auf die Seminarkonzeption, den Erwartungshorizont, und die zugrundeliegende Theorie zum Meaningful Learning nach Ausubel und den zugehörigen Forschungsergebnissen zu Concept Maps (vgl. Abschnitt 1.2). Von den vier Punkten des Erwartungshorizonts (vgl. Abbildung 4), der in Anlehnung an die Lernziele des Seminars erstellt wurde, können wir in dieser Analyse die Punkte 1 und 2 berücksichtigen. Die beiden jeweils zugehörigen Codes sind *gesellschaftliche Entwicklung* und *mathematische Denkweisen*. Der dritte Punkt des Erwartungshorizontes kann mit den Maps nicht erfasst werden und der vierte Punkt wurde hauptsächlich im nicht untersuchten Prüfungsteil abgefragt. Aus unseren theoriebezogenen Überlegungen wurden die Codes *node-link Diagramm*, *Liste/Tabelle* und *beschriftete Linie/Pfeil* entwickelt. Hierbei orientierten wir uns auch daran, inwiefern die von den Studierenden erstellten Knowledge-Maps sich von den gestalterisch eingeschränkten und mit rigideren Vorgaben versehenen Concept Maps unterscheiden.

In der inhaltlichen Dimension wurden insgesamt folgende Codes vergeben: *Fortschritte in der Mathematik* (A1), *Orientierung an den Seminarthemen* (A2), *gesellschaftliche Entwicklung* (A3), *mathematische Denkweisen* (A4), *Personen* (A5), *Zeit* (A6) und *Ort* (A7). In der gestalterischen Dimension wurden die Codes *Farben* (B1), *Symbole* (B2), *verschiedene Linienarten* (B3), *verschiedene Formen* (B4), *node-link Diagramm* (B5), *Zeitstrahl* (B6), *Liste/Tabelle* (B7), *geographische Karte* (B8) und *beschriftete Pfeile/Linien* (B9) vergeben.

In **Tabelle 4** (siehe Anhang) sind die Ergebnisse der Codierung der Knowledge-Maps hinsichtlich der auftretenden Elemente dargestellt. In der ersten Spalte befinden sich die Bezeichnungen der abgegebenen Knowledge-Maps, die mit anonymisierten Kürzeln der/des jeweiligen Studierenden versehen wurden. In der ersten Zeile stehen die Kürzel der Codes. Für die Codes bedeutet eine 1, dass die jeweilige Knowledge-Map entsprechend codiert wurde, eine 0 bedeutet, der Code wurde nicht vergeben. Die Tabelle ist so organisiert, dass eine Leerspalte die inhaltlichen von den gestalterischen Codes trennt. Leerzeilen wurden auch jeweils zwischen Maps, die in Form von Listen/Tabellen (B7) organisiert waren, einer Mischung aus Listen/Tabellen und node-link Diagrammen (B5 und B7) darstellten und node-link Diagrammen eingefügt.

Die **Antworten der Studierenden auf die erste Prüfungsfrage** („Was haben Sie sich bei der Erstellung der Knowledge-Map gedacht, wie sind Sie vorgegangen?") wurden deduktiv und induktiv codiert. Dabei stand die Frage, was die Studierenden über die Strukturierung der Knowledge-Map und ihr Vorgehen in Hinblick auf die Seminarziele bzw. den Erwartungshorizont sagten, im Vordergrund. Aus der Analyse der Knowledge-Maps wurden die Codes der inhaltlichen Dimension *Fortschritte in der Mathe-* *matik* (C1), *Orientierung an den Seminarthemen* (C2), *gesellschaftliche Entwicklung* (C3), *mathematische Denkweisen* (C4), *Personen* (C5), *Zeit* (C6) und *Ort* (C7) übernommen. Darüber hinaus sprachen die Studierenden neben gestalterischen Aspekten zusätzlich über die Technik des Zusammenfassens in Stichpunkten, induktiv codiert mit *Stichpunkte* (C8), und auftretende Schwierigkeiten bei der Maperstellung, und teilweise auch deren Lösungen, induktiv codiert mit *Schwierigkeiten* (C9).

Die Ergebnisse der Codierung wurden in **Tabelle 6** (siehe Anhang) zusammengefasst. Genau wie in Tabelle 4, wurden auch hier Leerzeilen eingefügt, um die unterschiedlichen Knowledge-Map-Formen zu kennzeichnen (kein node-link Diagramm, Mischung, node-link Diagramm). Die angesprochenen Schwierigkeiten und zugehörige Lösungen der Studierenden sind in **Tabelle 7** (siehe Anhang) zusammengefasst.

Da das vorliegende Projekt sich thematisch in den Bereichen Mathematikgeschichte und Mathematikdidaktik bewegt, waren wir inhaltlich in besonderem Maße an dem Potenzial interessiert, das der **Erarbeitung mathematischer Denkweisen zur Förderung von Reflexion** innewohnt. Aus diesem Grunde setzten wir Analyseergebnisse aus der Codierung der Knowledge-Maps mit den Analysen der Codierung der Transkripte der mündlichen Prüfung in Beziehung. Die Vermutung, dass das Erstellen der Knowledge-Map insbesondere in Hinblick auf die Entwicklung mathematischer Denkweisen Reflexion anregte, wird durch folgendes Zitat aus einer der Prüfungen unterstützt:

> HI: Und habe dann versucht, äh mit rot noch eine Abstraktionsebene weiter zu gehen und damit die mathematischen Denkweisen ähm zu kennzeichnen. Also was über diesen mathematischen Entwicklungen noch ähm drüber steht.

In der **Tabelle 5** (siehe Anhang) sind diese Ergebnisse zusammen mit den Protokollierungen der Interaktion der Studierenden mit der Knowledge-Map zusammengetragen. Die Einordnung der Maps in die Tabelle basiert auf der Codierung der Knowledge-Maps. Dabei wurden Maps, die node-link Diagramme darstellen (Code B5), die keine node-link Diagramme darstellen (Code B7) und Maps, die eine Mischung aus node-link Diagrammen und Listen darstellen (Code B5 und B7) unterschieden. Die Einordnung in die jeweiligen Zeilen basiert auf der Strukturierung der Knowledge-Map hinsichtlich mathematischer Denkweisen (Code A4), bzw. dem Fehlen einer solchen Strukturierung. Die hinter dem Kürzel angegebenen Buchstabenkombinationen stehen für die Codes aus den Transskripten: gesellschaftliche Entwicklung (C3), mathematische Denkweise (C4) und Bezug zu mathematik-historischen Personen (C5). In welchem Umfang die Studierenden mit ihrer eigenen Knowledge-Map während der mündlichen Prüfung interagierten (auf die Map deuten und auf die Map schauen), wird durch die Angabe der + und - Symbole angegeben. Dabei bedeutet ++ sehr viel, + viel und - wenig Interaktion der/des Studierenden mit der Map. Beispielsweise wurde die Knowledge-Map von AB mit kein node-link Diagramm (B7) und mathematische Denkweisen (A4) codiert. In der mündlichen Prüfung sprach diese Studentin dann im Rahmen der Antwort auf die erste Frage über mathematische Denkweisen (C4), gesellschaftliche Entwicklung (C3) und mathematik-historische Personen (C5). Während des gesamten ersten Prüfungsteils bezog sich die Studentin dabei viel (+) auf die Knowledge-Map.

Die **Fragebögen** wurden von allen 21 Studierenden, die an der mündlichen Prüfung teilnahmen, ausgefüllt und anschließend von uns elementar quantitativ ausgewertet. Illustrationen zu den Ergebnissen sind im Folgeabschnitt in den Text eingearbeitet.

### 3.3 Analyseergebnisse

*Codierung der Knowledge-Maps / Tabelle 4 im Anhang*

Insgesamt weisen die Knowledge-Maps eine große Vielfalt auf. Fünf der 21 Maps wurden nicht als node-link Diagramme codiert. Diese Maps waren entweder als Tabellen (AB und BC) oder als Stichpunktlisten (CD, EF und DE) organisiert (Code B7) und enthielten alle einen Zeitstrahl (Code B6). Drei Knowledge-Maps stellten eine Mischung aus node-link Diagrammen und Stichpunktlisten dar (FG, HG und HI). Die übrigen 13 Knowledge-Maps waren node-link Diagramme (Code B5). Beschriftete Pfeile oder Linien, ein zentrales Element bei Concept Maps (Code B9), traten in 11 der 21 Maps auf, wobei in den fünf Maps, die Listen oder Tabellen darstellten, gar keine verknüpfungsanzeigenden Linien oder Pfeile auftraten. Obwohl es strukturell naheliegend ist, haben nur acht Studierende ihre Knowledge-Map anhand der Seminarthemen (Code A2) organisiert. Alle Knowledge-Maps sind jedoch entweder nach Zeit (Code A6), nach Ort (Code A7) oder nach beidem organisiert, was bei diesem Seminar durch die Strukturierung der Themen (vgl. Abschnitt 2.2) nahegelegt wurde. Die beiden Codes gesellschaftliche Entwicklung (A3) und mathematische Denkweisen (A4), die sich aus dem Erwartungshorizont und den Seminarzielen ergaben, finden sich jeweils neun- bzw. 11-mal. Dabei kommen Knowledge-Maps, die mindestens eines der beiden Merkmale aufweisen, im Bereich kein node-link Diagramm, Mischung und node-link Diagramm vor (vgl. auch Tabelle 5). Bei den Codes aus der gestalterischen Dimension fällt auf, dass die Knowledge-Maps, die keine node-link

Diagramme darstellen weniger Designelemente aufweisen als node-link Diagramme oder Mischungen. Auch die eingesetzten Mittel zur Erstellung der Knowledge-Maps waren sehr vielfältig. 12 der abgegebenen Maps wurden digital, mit unterschiedlichen Programmen (neben verschiedenen Mind- und Concept-Mapping Programmen kamen auch Textverarbeitungs-, Datenbank- und Tabellen-Programme zum Einsatz) erstellt, die restlichen neun Knowledge-Maps wurden per Hand gezeichnet bzw. gemalt. Vor allem die per Hand erstellten Knowledge-Maps weisen zum Teil vielfältige kreative Aspekte auf, die durch die Codierung nicht erfasst werden.

*Ergebnisse der Codierung der Prüfungsausschnitte / Tabelle 6 im Anhang*

Betrachtet man die Codes C1 bis C7, die den Codes der inhaltlichen Dimension A1 bis A7 aus der Analyse der Knowledge-Maps entsprechen, so fällt auf, dass die Studierenden diese Aspekte in sehr unterschiedlichem Ausmaß angesprochen haben. Fünf Studierende haben nahezu alle Aspekte angesprochen (AB, HI, LM, ST und TU). Drei Studierende haben ausschließlich angesprochen, dass sie sich an den Seminarthemen orientiert haben (FG, RS und UV). Hinsichtlich der angesprochenen Aspekte gibt es keine Knowledge-Map Form, die hervorsticht, in allen drei Formen gab es jeweils Studierende die viele oder wenige Aspekte angesprochen haben. Fünf Studierende haben berichtet, dass sie zuerst Stichpunkte zusammengefasst haben und 14 Studierende berichteten von Schwierigkeiten bei der Erstellung der Map.

*Ergebnisse bezüglich der genannten Schwierigkeiten / Tabelle 7 im Anhang*

Eine zentrale Schwierigkeit betraf die Stofffülle und die daraus resultierende Unübersichtlichkeit der Verknüpfungen. Im Folgenden geben wir einige bereinigte Zitate der Studierenden aus dem ersten Teil der mündlichen Prüfung wieder.

> FG: […] am Anfang war das alles noch ganz leicht und ich konnte das auch gut alles nachschreiben. Aber dann am Ende wurden es immer mehr Leute und immer mehr Themen. […] die habe ich dann hier versucht draufzubringen.

Die Knowledge-Map von FG ist in besonderem Maße durch diese Schwierigkeit (immer mehr Leute und immer mehr Themen) und ihr Vorgehen (versucht das draufzubringen) bestimmt. Eine Lösung für das Problem hat sie nicht erarbeitet. Sie hat zu jedem Seminarthema eine Stichpunktliste erstellt und diese Listen als nodes in dem node-link Diagramm positioniert. Inhaltlich ist die Studentin nicht über die stichpunktartige Zusammenfassung hinausgekommen. Weder mathematische Denkweisen noch gesellschaftliche Entwicklungen sind in ihrer Knowledge-Map rekonstruierbar. Auch ihre Antwort auf die erste Prüfungsfrage offenbart keinerlei Gedanken bezüglich dieser beiden, oder sonstiger verbindender Elemente. Dass Schwierigkeiten mit dem Umfang und die Entscheidung gegen die Form node-link Diagramm aber nicht notwendig dazu führen, dass die Map inhaltlich oberflächlich bleibt, zeigt ein weiteres Beispiel:

> AB: Ja. Als erstes habe ich eine Mind Map gemacht. Und habe zu jedem Thema, was wir hatten, oder zu jeden Sitzungen einmal so ein, ja, so einen, so einen Ast gemacht. Und äh da dachte ich, gut, okay, äh da hm sieht man nicht so die Verknüpfung.

Die Studentin AB hat sich schließlich für ein Tabellendesign entschieden in dem sowohl gesellschaftliche Entwicklung als auch mathematische Denkweisen codiert wurden. Dabei hat sie diese (und andere) Aspekte explizit reflektierend erarbeitet:

> AB: Und denn habe ich gedacht, okay, was ist jetzt so das Wesentliche für mich?

Das mit dem Stoffumfang zusammenhängende Problem, dass womöglich der zur Verfügung stehende Platz nicht ausreicht, wurde ebenfalls von Studierenden angesprochen:

> QR: Ähm also das habe ich an so einer Whiteboard-Tafel erstmal gemacht, weil ich ja nicht sicher war, ob das dann so funktioniert.

Ebenso wurden Gedanken geäußert, bezüglich der grundlegenden Struktur, in welcher die Knowledge-Map am besten zu organisieren sei:

> MN: Ja. Also ich hatte mir zuerst überlegt, ob ich einen Zeitstrahl mache? Aber das ist halt, gerade weil Mesopotamien und Islam relativ viel überdecken, habe ich das in Länder organisiert.

> HI: Also in erster Linie halt überlegt, ob ich eben, also wichtige Kernprobleme der Mathematik raussuche und es danach ordne. Und habe mich dann aber doch entschieden, das Ganze chronologisch zu machen.

Außerdem wurde das Herausarbeiten der mathematischen Denkweisen als schwierig empfunden:

> NO: […] Ich habe versucht, so ein bisschen Denkweisen nahezubringen oder so vom Gefühl das Wichtigste der Denkweise. Aber das ist so ein bisschen schwierig, das zu kategorisieren und ähm deswegen habe ich auch oben in die Themen noch ein paar Pfeile reingebaut. Ähm und versucht, das so zu erklären.

*Ergebnisse bezüglich der Erarbeitung mathematischer Denkweisen zur Förderung von Reflexion / Tabelle 5 im Anhang*

An Tabelle 5 wird deutlich, dass Studierende, die eine Knowledge-Map erstellt hatten, welcher der Code mathematische Denkweise nicht zugeordnet werden konnte, auch in der mündlichen Prüfung nicht davon gesprochen haben, dass dies ein relevanter Teil ihres Vorgehens bei der Erstellung der Map war. Umgekehrt sprachen alle Studierende, deren Knowledge-Maps entsprechend codiert wurden auch in der mündlichen Prüfung über mathematische Denkweisen oder zumindest über gesellschaftliche Entwicklungen (BC und KL) als relevante Strukturierungsmerkmale. Für alle drei Map Formen gibt es sowohl Maps, denen mathematische Denkweise zugeordnet werden konnten, als auch Maps, bei denen das nicht möglich war. Insbesondere gibt es Maps, die keine node-link Diagramme darstellen, aber dennoch auf Reflexionsmomente, wie Herausarbeitung mathematische Denkweisen und gesellschaftlicher Entwicklungen, hindeuten, sowie Maps, die node-link Diagramme darstellen, die das nicht tun.

Bezüglich der Interaktion mit der Knowledge-Map in der mündlichen Prüfung kann festgestellt werden, dass bis auf CD nur Studierende, die node-Link Diagramme erstellt haben, sehr viel (++) mit der Map interagierten. Unter den Studierenden, die node-link Diagramme erstellt haben, befinden sich aber auch drei (TU, NO und RS), die nur wenig interagierten.

*Ergebnisse der Fragebögen*

Die Auswertung der Fragebögen ergab nachfolgende Ergebnisse hinsichtlich des Zeitpunkts der Erstellung der Knowledge-Maps durch die Studierenden: Obwohl von Semesterbeginn an thematisiert wurde, dass die Map als Prüfungsbegleitung der Memorierung aller Inhalte sowie der kognitiven Verknüpfung der einzelnen Referatsthemen dienen sollte, wurden lediglich drei Maps begleitend zum Seminar erstellt. Sechs von 21 Studierenden hatten die Knowledge-Map teilweise während des Verlaufs des Seminars konzipiert und vor der Prüfung fertiggestellt. Über die Hälfte der Studierenden fertigten ihre Map erst am Semesterende vor der Prüfung an (vgl. Tabelle 1).

| Parallel zum Seminar | Mischform | Am Ende des Semesters |
|---|---|---|
| 3/21 | 6/21 | 12/21 |

Tab. 1: Zeitpunkt der Erstellung der Knowledge-Map

Bei der zweiten Frage hinsichtlich des zeitlichen Aufwands variierten die Antworten der Studierenden stark: Neben fünf Studierenden, die eine qualitative (angemessen, hoch, sehr hoch) Bewertung abgaben, antworteten sieben Studierende mit Angaben zwischen drei und sechs Stunden, drei Studierende schätzen die Erstellung mit zehn bis zwölf Stunden ein und sechs Studierende sogar mit über zwanzig Stunden (vgl. Tabelle 2). Diese Antworten verdeutlichen, dass die Intensität der Auseinandersetzung mit der eigenen Map individuell stark unterschiedlich war.

| 3-6 h | 10-12 h | 20-25 h | Qualitative Angaben |
|---|---|---|---|
| 7/21 | 3/21 | 6/21 | 5/21 |

Tab. 2: Angegebener Arbeitsaufwand

Die dritte Frage nach der empfundenen Sinnhaftigkeit der Map-Erstellung parallel zum Seminarwurde sehr unterschiedlich beantwortet. Von 17 Studierenden wurde diese positiv eingeschätzt, von sieben Studierenden negativ (vgl. Tabelle 3).

| Maperstellung sinnvoll | Maperstellung nicht sinnvoll |
|---|---|
| 17/21 | 7/21 |

Tab. 3: Einschätzungen der Sinnhaftigkeit der Maperstellung parallel zum Seminar

Interessant war für uns, dass die negativen Antworten durchweg sehr ähnlich begründet wurden (vgl. Abb. 8): Die Studierenden empfanden es als schwierig, einen Gesamtüberblick über die Inhalte des Seminars aufzubauen bzw. direkt Verbindungen zwischen den einzelnen Seminarthemen herzustellen. Demgegenüber stehen sehr differenzierte und oft ausführliche Antworten bei einer positiven Einschätzung.

> Nein, weil einem dann der Gesamtüberblick und das „Ziel" des Seminars fehlt.

> Nein, da man erst später viele Verbindungen sieht.

> Ja, da man sich beim Erstellen sehr intensiv mit dem Thema beschäftigt und seinen eigenen Lernstil verfolgen kann wegen der wenigen Vorgaben.

> Ja, das festigt die Inhalte noch einmal. In den Vorträgen hat man viel Input auf einmal & kann die Verknüpfungen oft nicht direkt ziehen. Erst bei der Erstellung der Map ist mir ersichtlich geworden wie die verschiedenen Themenkreise & auch verschiedene Zeitintervalle zusammenhängen.

Abb. 8: Zwei negative und zwei positive Antworten auf die Frage zur Einschätzung der Sinnhaftigkeit

Explizit wurde die Freiheit bei der gestalterischen Umsetzung hervorgehoben und Merkmale des Meaningful Learning angesprochen: intensive Beschäftigung, eigener Lernstil, Festigung von Inhalten, Verknüpfung von inhaltlichen Themenkreisen und zeitlichen Abfolgen sowie die Erarbeitung eines Überblicks.

# 4. Diskussion der Ergebnisse

Im Folgenden sollen die Analyseergebnisse weiterführend diskutiert werden. Dieser Abschnitt ist, ähnlich wie auch Abschnitt 1, zweigeteilt. Im ersten Teil soll auf Basis eines mathematik-historischen Literaturbezugs die Dozentinnenperspektive auf den Einsatz unseres Knowledge-Map-Arbeitsauftrages dargelegt werden. Im Anschluss folgt die abschließende Bezugnahme auf die im ersten Abschnitt auf Basis der didaktischen Literatur aufgeworfenen Forschungsfragen.

## 4.1 Diskussion der Ergebnisse hinsichtlich mathematik-historischer Bezüge

Der Mathematikdidaktiker und -historiker Gregor Nickel (Universität Siegen) diskutiert in seinem Aufsatz „Vom Nutzen und Nachteil der Mathematikgeschichte für das Lehramtsstudium" (Nickel, 2013) positive und negative Aspekte beim Einsatz und der Vermittlung von mathematik-historischen Inhalten als hochschuldidaktisches Hilfsmittel. Hierbei fokussiert der Autor unterschiedliche Kategorien, die wir an dieser Stelle in Kürze anführen möchten, um uns in unserer Diskussion strukturell daran zu orientieren.

Als hilfreich bei der Erarbeitung eines besseren Verständnisses für Mathematik bewertet Nickel insbesondere vier Einsatzmöglichkeiten: den tröstenden, anekdotischen Gebrauch von Mathematikgeschichte, den exemplarischen Gebrauch, welcher als eine „Erfahrung mit echter Mathematik" (Nickel, 2013, S. 259) beschrieben wird, den Einsatz im (historisch-)genetischen Sinne (siehe u.a. Wagenschein, 1965, oder Kronfeller, 1998), sowie den verfremdenden Gebrauch, welcher von Nickel als „spezifisch auf das spätere Berufsfeld Schule" (Nickel, 2013, S. 258) bezogen betrachtet wird. Zu dem letzten Punkt sei erwähnt, dass der historische Kontext von Mathematik insofern als Verfremdung betrachtet werden kann, als dass Studierende zu vor Jahren selbst behandeltem Schulstoff nochmals die Erfahrung des eigenen Erlernens machen können. Diese alternative Perspektive kann insbesondere aufzeigen, „dass es keineswegs selbstverständlich und einfach ist, dass ‚man' so notiert und rechnet, wie es heute üblich ist." (Nickel, 2013, S. 259). Dem erlernten Blick auf mathematische Sachverhalte wird durch historische Inhalte im Optimalfall folglich eine zusätzliche reflektierende Ebene hinzugefügt.

Zwei weitere Varianten des Einsatzes charakterisiert Nickel als dem historischen Verständnis von Mathematik abträglich. Zum einen der Gebrauch in Form einer Beschreibung idealisierter Typisierungen, welcher als antiquarische Karikatur bezeichnet wird. Hiermit ist etwa die Aufzählung von verstaubten Kuriositäten der Mathematikgeschichte gemeint. Zum anderen der Gebrauch im monumentalistischen und jovialen Sinne, bei dem etwa unerreichbare Mathematikgenies in den Vordergrund gestellt werden. Darüber hinaus kann der Einsatz von Geschichte als Verständnishindernis fungieren, indem er beispielsweise bei Studierenden zu Verwirrung und zusätzliche Schwierigkeit zum regulären Stoff wahrgenommen wird.

Wir stellen uns die Frage, inwieweit die positiven Aspekte in unserer Studie durch die Einbindung von Knowledge Maps unterstützt wurden bzw. ob den negativen Aspekten entgegengewirkt werden konnte. Dafür sammeln wir im Folgenden Beobachtungen aus Dozentinnen- und Prüferinnen-Perspektive. Obwohl wir uns um eine möglichst objektive Stellungnahme bemühen, müssen diese, zumindest in Teilen, als subjektiv verstanden werden.

Hinsichtlich des **anekdotisch, tröstenden Gebrauchs** hebt Gregor Nickel unter anderem die biographische Einbeziehung von Persönlichkeiten aus der Mathematikgeschichte hervor, durch welche Studierenden unter anderem nahegebracht werden kann, dass stets ein langjähriger Entwicklungsprozess und ausdauernde Arbeit für mathematischen Fortschritt von Nöten waren. Dies kann als aufbauend bzw. tröstend wahrgenommen werden und die Frustration hinsichtlich eigener Verständnisprobleme senken. Biographien spielten auch bei den Referaten des Seminars eine wiederkehrende Rolle, wobei auffällig war, dass die Studierenden Lebenswege sehr unterschiedlich integrierten: Während einige Vortragende ausgewählte Biographien weitgehend getrennt von den mathematischen Inhalten an den Anfang ihres Seminarthemas stellten (beispielsweise bei Adam Ries), brachten andere die sozialen, gesellschaftlichen und politischen Lebensumstände von Persönlichkeiten mit ihrem Einfluss auf Mathematik in Verbindung (beispielsweise bei Sophie Germain mit starkem Fokus auf den Kontext der Französischen Revolution). Betrachten wir demgegenüber die erstellten Knowledge Maps (vgl. Codierung in Tabelle 4), so fällt auf, dass bei nahezu allen Maps (19/21) Persönlichkeiten (Code A5) als strukturierendes Element aufgenommen wurden und diese darüber hinaus mit anderen Kategorien verknüpft wurden. Lediglich in zwei von 21 Fällen (Knowledge-Maps LM und DE) stehen Namen isoliert von mathematik-historischen Ereignissen, wobei die Map DE zusätzlich zu den Namen noch Lebensdaten auflistet. Dies deckt sich mit unseren Beobachtungen aus den mündlichen Prüfungen. Die Erstellung der Knowledge Maps führte dazu, dass, unabhängig vom zugrundeliegenden Referat, die meisten Studierenden nach Verbindungen zwischen den mathematik-historischen Persönlichkeiten und anderen Aspekten such-ten und sich somit ein mehrschichtiges Bild von Mathematikgeschichte ergeben konnte. Auf die erste Frage in der mündlichen Prüfung gaben 12 von 21 Studierenden einen Bezug zu Personen an, davon auch einer der beiden Studierenden (LM, vgl. Tabelle 5), die Personen nicht als strukturierendes Element für die Knowledge-Map genutzt hatten. Unser Eindruck war insgesamt, dass insbesondere im Rahmen des Seminars der angesprochenen **Tendenz zur Monumentalisierung oder Heroisierung** von Persönlichkeiten und Entwicklungsschritten entgegengewirkt wurde.

Der **exemplarische Gebrauch** erlaubt nach Nickel „eine ‚Erfahrung Mathematik' im – zwar nicht aktuellen, aber doch authentischen – Forschungskontext" (Nickel, 2013, S. 259). Dazu, inwiefern die von den Studierenden hergestellten Forschungsbezüge als authentisch erfahren wurden haben wir keine Daten. Fortschritte in der Mathematik (Code A1) wurde in 14 von 21 Knowledge-Maps kodiert. Allerdings weist lediglich die Knowledge-Map MN eine auf Mathematik bezogene node („Zahlsysteme/neue Mathematik") auf, von der Verknüpfungen zu anderen nodes gehen. Die beiden Knowledge-Maps AB und BC sind als Tabelle organisiert, in der „Errungenschaften" jeweils eine Zeile kennzeichnen. In allen anderen Knowledge-Maps stehen aber Aspekte, die mit einem exemplarischen Gebrauch von Mathematikgeschichte in Verbindung gebracht werden können, eher im Hintergrund.

Wir stimmen Gregor Nickel zu, dass die Kenntnis der **historischen Genese mathematischer Denkweisen** ein „extrem anspruchsvolles Ziel" darstellt. Die Erarbeitung abstrakter Begriffe aus einer konkreten, phänomenbasierten Entwicklungsgeschichte heraus erfordert eine ungewohnte und mehrschichtige kognitive Auseinandersetzung. Im Rahmen der mündlichen Prüfungen haben wir die Beobachtung gemacht, dass diese schwierige Reflexion durch die Erstellung der Knowledge Maps zumindest unterstützt werden kann.

Als **Beispiel** kann hier **ST** angeführt werden: Die Studentin hat auf ihrer Knowledge-Map drei Verbindungsarten – Verbindungslinien, ein gezeichnetes „Zeitloch" und eine farbliche Kodierung – zu einem Themenkomplex miteinander kombiniert. Inhaltlich ging es dabei auf den ersten Blick um das so genannte Archimedes-Palimpsest, welches 1906 wiederentdeckt wurde und unter anderem „Die Methode", eine frühe Beschreibung eines heuristischen Approximationsverfahrens nach Archimedes, welches als wegweisend für die spätere Integralrechnung verstanden werden kann, enthält. Im Seminar wurde die Entdeckung und Entzifferung des historischen Fundes besprochen sowie auf die örtliche und zeitliche Verbreitung von Approximationsverfahren eingegangen. Insbesondere wurden zwei Verfahren zur Annäherung an eine Kreisfläche bzw. an die Kreiszahl Pi von Archimedes und Liu Hui gegenübergestellt. Die Strukturierung, welche von ST entwickelt wurde, entpuppte sich als erstaunlich komplex: Durch die Verwendung von gleichen Farben visualisiert die Studentin zunächst den „Sprung" zwischen der griechischen Antike (Archimedes) und dem 20. Jahrhundert, durch Verbindungslinien zwischen griechischer und chinesischer Antike stellt sie zudem die Übermittlung von Erkenntnissen dar und schließlich verknüpft sie verschiedene Flächenberechnungs- und Approximationsmethoden (Streifenmethode und Fläche unter einer Parabel, Approximation an Pi und Kreisinhalt) miteinander. Die überlegte Herangehensweise der Studentin an einzelne Aspekte der sehr umfassenden und vielseitigen Einflussgeschichte auf die Integralrechnung verdeutlicht die Möglichkeit der Sensibilisierung durch die Erstellung einer Knowledge-Map für eine historische Genese. Selbstverständlich ist hier kein Automatismus, sondern lediglich ein Angebot an die Studierenden gegeben.

Darüber hinaus wird gerade der vierte Punkt der **Verfremdung** und des daraus entstehenden möglichen neuen Prozesses des eigenen Lernens durch Verknüpfung von zuvor bekanntem Stoff (beispielsweise der Kreiszahl Pi als Verhältnis von Durchmesser und Radius eines Kreises) mit einer neuen Perspektive (beispielweise den beiden Approximationsverfahren an Pi durch Archimedes und Liu Hui) von der Methode der Knowledge-Maps angesprochen. Durch die bewusste eigene visuelle Vernetzung der Studierenden von größtenteils bereits bekannter (weitestgehend elementarer) Mathematik und in weiten Teilen bis dato unbekannten historischen Entwicklungsschritten war eine Vertiefung des inhaltlichen Verständnisses bemerkbar. Dies wurde im ersten Teil des Prüfungsgesprächs an wiederholter Stelle deutlich, was wir anhand von zwei Beispielen veranschaulichen möchten:

Die **Studentin JK** ordnete ihre Knowledge-Map als „Stammbaum der Mathematik" (eigene Bezeichnung der Studentin) an, bei dem sie mehrfach dezidiert astronomische Erkenntnisse hervorhob. Dies begründete sie in der mündlichen Prüfung damit, dass auch nachdem man historisch betrachtet,

> nicht mehr nur auf praktisches Rechnen oder Mathematik als Hilfsmittel angewiesen war […] zum Beispiel noch bei Gauß [ersichtlich wurde] dass er dort Asteroidenbahnen berechnen wollte.

Sie berücksichtigte diese Erkenntnis, indem sie eine Unterkategorie „Astronomie" aufmachte und diese nutzte, um ihr bisheriges Verständnis von der Unterteilung in eine antike, praxisbezogene Mathematik und eine neuere, theoriebezogene Mathematik aufzubrechen.

Der **Student EF** legte einen inhaltlichen Schwerpunkt auf diophantische Gleichungen und die Fermatsche Vermutung. Die Anordnung seiner Knowledge-Map in Stichpunktlisten zu den einzelnen Seminarthemen sowie einem durchgängigen Zeitstrahl anhand der Lebensdaten von einzelnen MathematikerInnen erschien auf den ersten Blick wenig vernetzend bzw. komplex. Auf Andrew Wiles angesprochen, kommentierte der Student allerdings seine Map dahingehend, dass der

> gesamte untere Zeitstrahl mit [hier] Andrew Wiles, den müsste man ja einen ganz dicken Strich […] zu Pierre de Fermat machen und letztendlich dann ja auch wieder zurück zur griechischen Antike. […] ziehen. Er betonte folglich, mit Verweis auf seine eigens angefertigte Knowledge-Map bzw. darauf, dass er diese entscheidende Verknüpfung nicht eingezeichnet hatte, seine Auffassung, dass inhaltlich ein breiter historischer Bogen zwischen modernen mathematischen Erkenntnissen (Beweis der Fermatschen Vermutung von 1994) und antiker Mathematik (diophantische Gleichungen) gezogen werden kann. An dieser Stelle sei erwähnt, dass im Seminar deutlich wurde, dass der Beweis von Wiles, welcher auf der Verknüpfung von mathematischen Forschungsgebieten (u.a. der Theorie der elliptischen Kurven) beruht, in seiner Komplexität das Wissen der Lehramtsstudierenden überstieg. Der Rückbezug zu elementaren diophantischen Gleichungen stellte für den Studenten anscheinend eine neue ermutigende Erkenntnis dar. Interessant ist, dass der Student, obwohl er Schwierigkeiten mit einer vernetzten Anordnung der Inhalte hatte (vgl. auch Tabelle 7), durch die Erstellung der Knowledge-Map eine vernetzende Auseinandersetzung mit den Inhalten erkennen lässt.

Ähnliche Beispiele, welche verdeutlichen, dass die Knowledge-Maps durch die Auseinandersetzung mit visuellen Kennzeichen von inhaltlichen und zeitlichen Verknüpfungen (auch, wenn festgestellt wurde, dass dies gerade nicht gelang!) das eigene Verständnis von Zusammenhängen unterstützen, gab es in nahezu allen mündlichen Prüfungen. Wir sehen folglich den positiven Aspekt, dass insbesondere in Teilen bekannte mathematische Sachverhalte elaboriert und neu betrachtet wurden.

Anhand der aufgeführten Beispiele und unseren Beobachtungen im Seminar sowie bei den mündlichen Prüfungen möchten wir insgesamt festhalten, dass wir die Knowledge-Maps als eindeutig unterstützend hinsichtlich der von Gregor Nickel in den Fokus genommenen Aspekte bewerten. Natürlich liefern diese keine Garantie hinsichtlich einer sensiblen transdisziplinären und strukturierten Auseinandersetzung mit Mathematikgeschichte, sie können diese allerdings merklich unterstützen.

Die von Nickel angesprochenen Verständnishindernisse beim Einsatz von Mathematikgeschichte konnten wir ebenfalls beobachten. Einige Studierende berichteten über unterschiedliche Schwierigkeiten bei der Erstellung der Knowledge-Maps (vgl. Tabelle 7). Dabei spielte auch die zum Teil komplexe Verknüpftheit und die Prozesshaftigkeit geschichtlicher Entwicklungen eine Rolle. Da wir diesen Punkt nicht systematisch erfasst haben, haben wir keinen vollständigen Überblick über alle aufgetretenen Schwierigkeiten. Die meisten der Studierenden, die über Probleme berichteten, entwickelten aber eigenständig Lösungen. Möglicherweise können diese Verständnishindernisse also als kognitive Hürden wirken, deren Überwindung potentiell Lern- und Reflexionsprozesse anregt.

## 4.2 Diskussion der Ergebnisse hinsichtlich didaktischer Bezüge

In diesem Abschnitt möchten wir unsere Ergebnisse in Bezug setzen zu den theoretischen Überlegungen aus Abschnitt 1.2.

*Zum Lernbegleitenden Einsatz der Knowledge-Maps*

Der Versuch des Einsatzes von Knowledge-Maps als prozessbegleitende Lernhilfe scheiterte in unserem Fall daran, dass eine der notwendigen Bedingungen für Meaningful Learning nicht erfüllt war, nämlich die Kooperationsbereitschaft der Studierenden. Trotz des klaren Arbeitsauftrages und mehrmaliger Empfehlungen erstellte die Mehrheit der Studierenden die Knowledge-Map erst gegen Ende des Semesters.

Die Ursachen hierfür könnten einerseits in den Rahmenbedingungen des Studiums begründet liegen, wie zum Beispiel Arbeitsgewohnheiten, die aus anderen Veranstaltungen übernommen werden, hohe Arbeitsbelastung während der Vorlesungszeit, welche einen Aufschub von Arbeit bis nach Ende der Vorlesungszeit nahelegt, etc.. Andererseits wurden wir von Studierenden auf Schwierigkeiten mit dem konkreten Arbeitsauftrag hingewiesen. Im Fragebogen und während der mündlichen Prüfungen wurde das Problem geäußert, dass man während des Semesters noch nicht genau wusste, welche Informationen über den Rest des Semesters noch hinzukommen würden. Gerade für diejenigen, die per Hand zeichneten, bedeutete dies eine große Unklarheit darüber, wie der Platz auf dem Blatt am besten einzuteilen sei, und wieviel Platz insgesamt benötigt würde. Dieses Problem hängt sicherlich auch mit der großen Fülle an Inhalten zusammen, die auf die Map gebracht werden sollten. Concept Maps sind primär dafür gedacht, einzelne Konzepte abzubilden, unsere Knowledge-Maps hingegen sollten die gesammelten Seminarinhalte enthalten.

Wir schlagen daher für den lernbegleitenden Einsatz von Knowledge-Maps vor, eine Einschränkung des Arbeitsauftrages zu erwägen. Sowohl eine Beschränkung auf eine überschaubare Menge an abzubildendem Stoff, als auch auf digitale Mapping-Werkzeuge[4] zur Erstellung der Map erscheinen hier zielführend. Gerade digitale Hilfsmittel eröffnen die Möglichkeit einer nachträglichen Umstrukturierung der Map, und könnten dadurch Reflexion unterstützen. Eine weitere Maßnahme, die hilfreich sein könnte, besteht in der Vorbereitung der Studierenden auf die Erstellung von Knowledge-Maps zu Beginn einer Veranstaltung, und der Erarbeitung von Strategien, mit dem Problem der Platzeinteilung und der Strukturierung der Map umzugehen.

*Maps zur Zusammenfassung und zum Lernen*

Der Einsatz der Knowledge Map zum Zusammenfassen von Informationen im Vorfeld einer Prüfung war die Art des Gebrauchs, die bei beinahe allen Studierenden zu beobachten war. Die Relevanz der Knowledge-Map für die Prüfung und der Nutzen der Maperstellung mit Blick auf die Prüfung, scheint für viele Studierenden eine gute Motivation gewesen zu sein, den Arbeitsauftrag gewissenhaft zu bearbeiten (z.B. Abbildung 9).

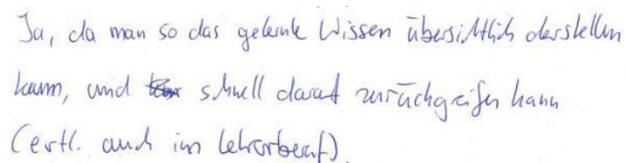

Abb. 9 Studierendenzitat zum Nutzen der Knowledge-Map zum Zusammenfassen und zum Lernen.

Obwohl, wie bereits erwähnt, im Rahmen unseres Arbeitsauftrages Knowledge-Maps entstanden, die sowohl inhaltlich als auch vom Format her Unterschiede zu den Concept Maps von Novak und Gowin aufweisen, sind wir der Ansicht, einige Aussagen aus entsprechender Forschungsliteratur auf unseren Arbeitsauftrag übertragen zu können.

So sind wir etwa der Ansicht, beobachtet zu haben, dass unser Arbeitsauftrag Reflexion über die Seminarinhalte ausgelöst hat. Das Auftreten von übergeordneten Konzepten wie mathematischen Denkweisen und gesellschaftlichen Entwicklungen in den Knowledge-Maps lässt vermuten, dass im Vorfeld bewusste Entscheidungen in Bezug auf die Art der Integration getroffen werden mussten. Entscheidungen, die zumindest einige Studierende erst einmal vor Probleme stellten (siehe Tabelle 7). Somit erscheint unser Arbeitsauftrag gerade aufgrund der Fülle an Informationen und deren Überkomplexität in der Lage, Problemlöseprozesse in Gang zu setzen, die nicht zuletzt die Struktur des behandelten Stoffes zum Inhalt haben. Auf der einen Seite bekräftigt dies die in 1.2.2 zitierte Ansicht von Novak und Gowin über die Wirkweise von Concept Maps:

> (…) the act of mapping is a creative activity, in which the learner must exert effort to clarify meanings, by identifying important concepts, relationships, and structure within a specified domain of knowledge. (Cañas et al., 2003, S. 22)

Auf der anderen Seite unterstützt diese Beobachtung vor dem Hintergrund fehlender Formatvorgaben in unserem Arbeitsauftrag und der Tatsache, dass übergeordnete Konzepte in allen Mapformaten auftreten konnten, auch die Vermutung von Nesbit und Adesope (2006), dass die genauen Formatvorgaben von Concept Maps möglicherweise nicht den entscheidenden Einfluss auf das Lernen ausmachen, sondern der generell höhere Grad an Aktivierung.

Was die im vorangegangenen Abschnitt empfohlene Modifikation des Arbeitsauftrages angeht, könnte man nun spekulieren, ob die für einen lernbegleitenden Einsatz von Knowledge-Maps vorgeschlagene Beschränkung studentischer Kreativität durch das Vorschreiben eines spezifischen Programms zur Erstellung der Map möglicherweise kontraproduktiv auf das Herbeiführen eines Problemlöseprozesses wirkt. Außerdem haben zumindest einige Studierende die „wenigen Vorgaben" ganz konkret positiv hervorgehoben (vgl. Abbildung 8). An dieser Stelle sei ferner angemerkt, dass sich unser Arbeitsauftrag trotz seines großen Umfangs bisher für den Einsatz in Universitätsseminaren als geeignet erwies. Wir empfehlen jedoch Vorsicht bei einer Übertragung in Kurse mit jüngeren Studierenden, die mit dem hiesigen Auftrag überfordert sein könnten.

Was die Akzeptanz des Arbeitsauftrages durch die Studierenden angeht, so gab es überwiegend positive Resonanz. Unter den Rückmeldungen waren einige, die die Methode als überaus hilfreich bewerteten, darunter zum Beispiel ein Student des Lehramts mit der Fächerkombination Mathematik und Philosophie, der uns nach der Prüfung mitteilte, dass er die Methode vorher nicht kannte, sie aber jetzt zum Lernen in Philosophieveranstaltungen nutzte. Vereinzelte Studierende standen dem Ganzen aber auch indifferent bis kritisch gegenüber (vgl. Abbildung 8). Aus unserer Sicht ist es gerade für Lehramtsstudierende wichtig, in ihrem Studium unterschiedliche Methoden kennenzulernen. Dabei halten wir es aus didaktischer Sicht für sinnvoll, die Studierenden die Methode nach Möglichkeit selbst ausprobieren zu lassen.

Durch unsere Art des Einsatzes der Maps kann dies im Rahmen eines nicht methodenfokussierten Seminars „einfach so nebenbei" erfolgen.

*Evaluation der Map: Bewertung und Diagnostik*

Ein diagnostisches Moment findet sich in unserem Prüfungsablauf an der Stelle, an der die Dozentin Fragen zur Knowledge-Map generiert.

Im Unterschied zu herkömmlichen Methoden, Fragen aus dem behandelten Stoff oder den Erwartungen der/des DozentInnen abzuleiten, knüpft unser Verfahren direkt an die Wissensstrukturen der Studierenden an. Speziell zur Entwicklung von Fragen in Bezug auf die Reflexion eigener Fehlvorstellungen (siehe Anforderungskategorie 3) entpuppten sich unsere Knowledge-Maps aus Sicht der Prüferinnen als äußerst hilfreich. Sie halfen während der Vorbereitung der Prüfung dabei, Verständnislücken bei den Studierenden zu erkennen, und die Prüfungsfragen daran zu orientieren. Gleichzeitig dienten die Vorstellungen der Studierenden, die sich in den Knowledge-Maps zeigten, der Dozentin als Anhaltspunkt für eine Reflexion ihrer eigenen Veranstaltung, der behandelten Inhalte und Vermittlungsmethoden. Das Auftreten der gleichen Fehlvorstellungen oder Wissenslücken in mehreren Maps kann in diesem Kontext als zielgebende Anregung verstanden werden, für eine mögliche Überarbeitung der Veranstaltung. Ein erwähnenswerter negativer Punkt dieses prüfungsvorbereitenden Einsatzes von Knowledge-Maps ist der relativ hohe Arbeitsaufwand, der im Vergleich zu herkömmlichen Arbeitsmodi entsteht. Wir ziehen einen solchen Einsatz von Knowledge-Maps in unserer Lehre jedoch weiterhin in Erwägung.

Im Vorfeld unseres Seminars hatten wir uns bewusst gegen die Bewertung der Knowledge-Maps selbst entschieden. Wie in Abschnitt 1.2 angedeutet, wurde in der Literatur Skepsis geäußert, dass eine möglichst objektive numerische Kategorisierung einzelner Concept Maps vor dem Hintergrund einer Bewertungssituation möglich bzw. sinnvoll wäre. Die Erfahrungen aus der Masterarbeit zur Veranstaltung in Wuppertal deuteten in die gleiche Richtung. Daher war unser Ansatz an dieser Stelle zurückhaltender und wir verwendeten die Knowledge-Maps lediglich als Leitfaden und Hilfsinstrument für die eigentliche Prüfung, sowie als unbenotete Studienleistung.

Eine Beobachtung, welche uns in dieser Zurückhaltung weiter bestärkt, ist die, dass wir anhand der Knowledge-Maps keine Vorhersagen über den Verlauf der mündlichen Prüfungen machen konnten, bzw. dass auch keine Korrelation zwischen den notenmäßigen Ergebnissen der Prüfungen und den Maps ins Auge gesprungen ist. Um diese ersten Beobachtungen genauer zu prüfen, wären aber größer angelegte Studien nötig, die außerdem sicherstellen müssten, dass die Ergebnisse einer mündlichen Prüfung hinreichend reliabel sind. Insofern können wir an dieser Stelle keine definitiven Aussagen treffen.

*Zum Einsatz der Map in der mündlichen Prüfung*

Im Vorfeld des Einsatzes in den Prüfungen vermuteten wir, dass es dem/der zu Prüfenden mit der Knowledge-Map leichter fallen könnte, in ein Gespräch einzusteigen. Außerdem erschienen Situationen, in denen Unklarheit bezüglich einer Frage entsteht und das Gespräch deswegen stagniert, bei diesem Format unwahrscheinlicher.

Tatsächlich stellte sich heraus, dass die Knowledge-Maps für die Dozentin nach entsprechender Vorbereitung eine Gesprächserleichterung boten. Insbesondere das Aufgreifen von Metaphern der Studierenden wurde als eine Kommunikationserleichterung empfunden. Die visuelle Unterstützung durch die Maps half nicht nur dabei, Verständnislücken bei den Studierenden zu erkennen, sondern auch, diese konkret zu benennen und mit den zu Prüfenden zu diskutieren. In solchen Diskussionen tritt die Arbeitsweise der Studierenden zutage und kann in der Bewertung berücksichtigt werden.

Die Studierenden nutzten ihre Knowledge-Maps während der Prüfung sehr unterschiedlich (vgl. Tabelle 5, Interaktion). Einige Studierende interagierten sehr häufig (durch verbales oder tatsächliches Deuten, durch Verweis auf spezifische Merkmale) mit ihrer Map, während andere Skrupel zu haben schienen, ihre Knowledge-Map allzu offensichtlich anzusehen. Hier hätten wir deutlicher Kommunizieren müssen, dass eine aktive Verwendung der Map in der Prüfung von uns nicht als „Spicken" verstanden wird, wie es herkömmliche Gewohnheiten in Leistungsüberprüfungen möglicherweise nahelegen. In unseren Daten wurden außerdem node-link Diagramme von den Studierenden häufiger als Gesprächsfokus genutzt, als andere Arten von Knowledge-Maps. Es stellt sich die Frage, ob dies rein zufällig geschah, oder ob das Map-Format einen Einfluss auf Interaktionen über die Map hat. Diese Frage kann auf Basis unserer Datenlage jedoch nicht weiter geklärt werden.

Was uns überraschte war, dass Studierende nicht immer zu eigenen Notizen Stellung nehmen konnten. Insofern scheint eine Gesprächserleichterung auf Studierendenseite durch den Einsatz von Knowledge-Maps nicht garantiert zu sein.

### 4.3 Fazit und Ausblick

In unserer Rolle als Dozentinnen werten wir den Einsatz unseres Knowledge-Map Auftrages als Erfolg. Wir konnten den Lehramtsstudierenden eine neue Methode beibringen, regten sie zur Auseinandersetzung mit mathematik-historischen Arbeits- und Denkweisen an, und konnten die Maps zusätzlich als Feedback für unser Veranstaltungsdesign nutzen. Wir zögern nicht, die in diesem Artikel beschriebene Methode auch weiterhin in unseren mathematik-geschichtlichen und mathematikdidaktischen Seminaren einzusetzen.

Was eine mögliche Modifikation des von uns verwendeten Arbeitsauftrages angeht, wollen wir keine pauschale Empfehlung aussprechen. Für einen lernbegleitenden Einsatz von Knowledge-Maps, oder für eine jüngere Zielgruppe, halten wir eine Einschränkung des Auftrages für sinnvoll. Zum Fördern von Problemlöseverhalten in Universitätsseminaren und um eine tiefergehende Auseinandersetzung mit dem Stoff anzuregen, scheint der aktuelle Auftrag gerade aufgrund seines hohen Anspruches in Besonderer Weise geeignet.

Was die Bewertung studentischer Maps angeht, so ist klar, dass in Veranstaltungen mit großer Teilnehmerzahl, wie etwa Vorlesungen, keine Einbindung von Knowledge-Maps in einen Benotungsprozess nach hiesiger Facon erfolgen kann. Es gibt aber universitäre Großveranstaltungen in den Bereichen Mathematikdidaktik und -geschichte, die eine abschließende Prüfung erfordern. Dies brachte uns auf den Gedanken, dass man große Zahlen von Knowledge-Maps eventuell mit Hilfe eines Comparative-Judgement-Verfahrens (etwa Bisson, Gilmore, Inglis & Jones, 2016) scoren könnte. Um dies überhaupt möglich zu machen, müssten jedoch einschneidende Einschränkungen des Arbeitsauftrages vorgenommen werden. Im Moment werden Abschlussarbeiten zu der Frage betreut/geschrieben, ob in diese Richtung Potenzial bestünde, oder nicht.

In unserer Rolle als Forschende stellten wir fest, dass sich einige Hypothesen aus der Literatur über die Wirkweise von Concept Maps in unseren qualitativen Daten zu bestätigen scheinen. Des Weiteren ergeben sich eine Reihe von Fragen zur weiteren Beforschung von Knowledge-Maps in unterschiedlichen Richtungen:

- Welche Arten von Problemen treten beim Erstellen der Knowledge-Maps (lernbegleitend oder zur Prüfungsvorbereitung) auf und inwiefern befördert das Lösen dieser Probleme, möglicherweise durch den Einsatz gezielter Unterstützungsmaßnahmen Seitens der Dozierenden, Lern- und Reflexionsprozesse?

- Hat das Mapformat (node-link Diagramm, kein node-link Diagramm, Mischung) Einfluss auf die Interaktion der Studierenden mit ihren Maps in der mündlichen Prüfung?

- Hat die Knowledge-Map Einfluss auf den Gesprächsverlauf in der mündlichen Prüfung? Zeigen sich insbesondere Auswirkungen auf den Umgang mit Themenwechsel?

## Anmerkungen

[1] Siehe beispielsweise die Projekte: Mathematics Genealogy Project und „Mapping the Republic of Letters" auf http://republicofletters.stanford.edu

[2] Nesbit und Adesope verstehen darunter „two-dimensional visual knowledge-representations, including flowcharts, timelines, and tables, that show relationships among concepts or processes by means of spatial position, connecting lines, and intersecting figures" (2006, S. 413).

[3] Für die Analyse haben wir das open access Online-Tool QCAmap, siehe auch www.qcamap.org, genutzt.

[4] Ein unserer Ansicht nach geeignetes Programm zur Erstellung von Concept Maps ist CmapTools, https://cmap.ihmc.us/, vgl. auch (Cañas et al., 2004).

**Anschrift der Verfasser**

Nicola Oswald
Leibniz Universität Hannover
Institut für Didaktik der Mathematik und Physik
AG Didaktik der Mathematik
Welfengarten 1
30167 Hannover
oswald@idmp.uni-hannovr.de

Sarah Khellaf
Leibniz Universität Hannover
Institut für Didaktik der Mathematik und Physik
AG Didaktik der Mathematik
Welfengarten 1
30167 Hannover
khllaf@idmp.uni-hannovr.de

Jana Peters
Leibniz Universität Hannover
Institut für Didaktik der Mathematik und Physik
AG Didaktik der Mathematik
Welfengarten 1
30167 Hannover
peters@idmp.uni-hannovr.de

**Anhang**

| Dokument | A1 | A2 | A3 | A4 | A5 | A6 | A7 | B1 | B2 | B3 | B4 | B5 | B6 | B7 | B8 | B9 |
|---|---|---|---|---|---|---|---|---|---|---|---|---|---|---|---|---|
| AB | 1 | 0 | 1 | 1 | 1 | 1 | 1 | 0 | 0 | 0 | 0 | 0 | 1 | 1 | 0 | 0 |
| BC | 1 | 1 | 1 | 1 | 1 | 1 | 0 | 0 | 1 | 0 | 0 | 0 | 1 | 1 | 0 | 0 |
| CD | 0 | 1 | 0 | 0 | 1 | 1 | 0 | 1 | 0 | 0 | 0 | 0 | 1 | 1 | 0 | 0 |
| DE | 0 | 0 | 0 | 0 | 0 | 1 | 1 | 1 | 0 | 0 | 0 | 0 | 1 | 1 | 1 | 0 |
| EF | 0 | 1 | 0 | 0 | 1 | 1 | 0 | 1 | 0 | 0 | 0 | 0 | 1 | 1 | 0 | 0 |
| FG | 0 | 1 | 0 | 0 | 1 | 1 | 0 | 1 | 0 | 0 | 0 | 1 | 0 | 1 | 0 | 1 |
| GH | 1 | 0 | 1 | 0 | 1 | 0 | 1 | 1 | 1 | 1 | 1 | 1 | 0 | 1 | 0 | 1 |
| HI | 1 | 0 | 1 | 1 | 1 | 1 | 0 | 1 | 0 | 0 | 1 | 1 | 1 | 1 | 0 | 1 |
| IJ | 0 | 0 | 0 | 0 | 1 | 1 | 0 | 1 | 0 | 1 | 1 | 1 | 1 | 0 | 0 | 0 |
| JK | 1 | 0 | 0 | 1 | 1 | 1 | 0 | 1 | 0 | 1 | 0 | 1 | 1 | 0 | 0 | 0 |
| KL | 1 | 0 | 1 | 1 | 1 | 1 | 0 | 1 | 1 | 1 | 1 | 1 | 1 | 0 | 0 | 1 |
| LM | 1 | 0 | 0 | 1 | 0 | 0 | 1 | 0 | 0 | 0 | 1 | 1 | 0 | 0 | 0 | 1 |
| MN | 1 | 0 | 0 | 0 | 1 | 0 | 1 | 1 | 0 | 0 | 0 | 1 | 0 | 0 | 0 | 1 |
| NO | 1 | 0 | 0 | 1 | 1 | 1 | 1 | 1 | 0 | 1 | 1 | 1 | 0 | 0 | 0 | 1 |
| OP | 1 | 1 | 1 | 1 | 1 | 1 | 0 | 1 | 1 | 0 | 0 | 1 | 1 | 0 | 0 | 1 |
| PQ | 1 | 0 | 0 | 0 | 1 | 1 | 1 | 1 | 0 | 1 | 1 | 1 | 0 | 0 | 0 | 0 |
| QR | 1 | 1 | 1 | 1 | 1 | 1 | 1 | 1 | 1 | 1 | 0 | 1 | 0 | 0 | 0 | 1 |
| RS | 0 | 1 | 0 | 0 | 1 | 1 | 0 | 1 | 0 | 1 | 0 | 1 | 0 | 0 | 0 | 1 |
| ST | 1 | 0 | 1 | 1 | 1 | 1 | 1 | 1 | 0 | 0 | 1 | 1 | 0 | 0 | 0 | 0 |
| TU | 0 | 0 | 1 | 1 | 1 | 1 | 0 | 1 | 0 | 0 | 1 | 1 | 1 | 0 | 0 | 0 |
| UV | 1 | 1 | 0 | 0 | 1 | 1 | 0 | 0 | 0 | 0 | 1 | 1 | 1 | 0 | 0 | 1 |

Tab. 4 In den Knowledge-Maps codierte Elemente

| | | kein node-link Diagramm | | | Mischung | | | node-link Diagramm | |
|---|---|---|---|---|---|---|---|---|---|
| | AB | (C4, C3, C5) | + | HI | (C4, C3, P) | - | JK | (C4, C3, C5) | ++ |
| | BC | (C4) | + | | | | ST | (C4, C3, C5) | ++ |
| | | | | | | | LM | (C3, C5) | ++ |
| **mathematische** | | | | | | | OP | (C3) | ++ |
| **Denkweise** | | | | | | | QR | (C3) | ++ |
| | | | | | | | KL | (C4) | ++ |
| | | | | | | | TU | (C4, C3, C5) | - |
| | | | | | | | NO | (C3, C5) | - |
| | CD | (C5) | ++ | GH | (C5) | + | PQ | (C5) | ++ |
| **keine** | EF | (C5) | + | FG | (-) | - | IJ | (C5) | ++ |
| **mathematische** | DE | (-) | - | | | | UV | (-) | ++ |
| **Denkweise** | | | | | | | MN | (-) | + |
| | | | | | | | RS | (-) | - |

Tab. 5 Kombination der Analyseergebnisse aus den Knowledge-Maps und den Prüfungsgesprächen

| Dokument | C1 | C2 | C3 | C4 | C5 | C6 | C7 | C8 | C9 |
|---|---|---|---|---|---|---|---|---|---|
| **AB** | 1 | 0 | 1 | 1 | 1 | 1 | 1 | 0 | 1 |
| **BC** | 1 | 0 | 1 | 0 | 0 | 1 | 0 | 0 | 1 |
| **CD** | 0 | 0 | 0 | 0 | 1 | 1 | 0 | 1 | 0 |
| **DE** | 0 | 0 | 0 | 0 | 0 | 1 | 1 | 0 | 0 |
| **EF** | 0 | 1 | 0 | 0 | 1 | 0 | 0 | 1 | 1 |
| **FG** | 0 | 1 | 0 | 0 | 0 | 0 | 0 | 1 | 1 |
| **GH** | 0 | 0 | 0 | 0 | 1 | 1 | 1 | 0 | 1 |
| **HI** | 1 | 0 | 1 | 1 | 1 | 1 | 0 | 0 | 1 |
| **IJ** | 1 | 0 | 0 | 0 | 1 | 1 | 0 | 0 | 1 |
| **JK** | 0 | 1 | 1 | 1 | 1 | 0 | 0 | 0 | 0 |
| **KL** | 1 | 1 | 1 | 0 | 0 | 0 | 0 | 0 | 0 |
| **LM** | 1 | 1 | 0 | 1 | 1 | 0 | 1 | 0 | 0 |
| **MN** | 1 | 0 | 0 | 0 | 0 | 0 | 1 | 0 | 1 |
| **NO** | 1 | 1 | 0 | 1 | 1 | 0 | 0 | 0 | 1 |
| **OP** | 1 | 1 | 0 | 1 | 0 | 0 | 0 | 0 | 1 |
| **PQ** | 0 | 0 | 0 | 0 | 1 | 1 | 1 | 0 | 0 |
| **QR** | 1 | 0 | 0 | 1 | 0 | 1 | 1 | 1 | 1 |
| **RS** | 0 | 1 | 0 | 0 | 0 | 0 | 0 | 0 | 1 |
| **ST** | 1 | 0 | 1 | 1 | 1 | 1 | 1 | 1 | 1 |
| **TU** | 1 | 0 | 1 | 1 | 1 | 1 | 0 | 0 | 1 |
| **UV** | 0 | 1 | 0 | 0 | 0 | 0 | 0 | 0 | 0 |

Tab. 6 In den transkribierten Prüfungsgesprächen codierte Elemente

| Schwierigkeiten | Lösung | |
|---|---|---|
| Stofffülle und Unübersichtlichkeit der Verknüpfungen | Weglassen (Tabelle/Liste) | EF, AB |
| | Farbe | IJ, RS, OP |
| | andere Map-Organisation/ Reduktion | GH |
| | Zoom | TU |
| | keine Lösung | FG |
| Unsicherheit bezüglich des zur Verfügung stehenden Platzes | Whiteboard | QR |
| | einfach versucht | NO |
| | keine Lösung | ST |
| Struktur der Map unklar | mehrere Modelle probiert | BC, MN, HI |
| Schwierigkeiten Mathematische Denkweisen zu charakterisieren | keine Lösung | NO |

Tab. 7 Von den Studierenden berichtete Schwierigkeiten bei der Erstellung der Knowledge-Map